\documentclass[preprint,10pt,3p]{elsarticlePLF}
\usepackage{amsfonts,amsmath}
\usepackage{tikz}
\usepackage{epsfig}
\usepackage{epstopdf}
\usetikzlibrary{shadings}
\usetikzlibrary{arrows}
\usetikzlibrary{shapes,snakes}
\usepackage{bm}
\usepackage{mathtools}
\usepackage{graphicx}
\usepackage{caption}
\usepackage[colorlinks=true ]{hyperref}
\usepackage{amssymb}

\usepackage[final]{pdfpages}

\newtheorem{definition}{Definition}[section]

\newtheorem{Claim}[definition]{Claim}
\newtheorem{remark}[definition]{Remark}
\newtheorem{example}[definition]{Example}
\newtheorem{proof}[definition]{Proof}

\include{TeXFigs}

\newcommand \bei {\begin{itemize}}
\newcommand \eei {\end{itemize}}
\newcommand \ubar u
\newcommand \RR {\mathbb R} 
\newcommand \del \partial

\newcommand \la \langle
\newcommand \ra \rangle 

\newcommand \auth    \textsc
\newcommand \be {\begin{equation}}
\newcommand \ee {\end{equation}}
\newcommand \bel {\be\label}
\newcommand \bcor {\begin{Claim}}
\newcommand \ecor {\end{Claim}}

\newcommand \bcorl {\bcor\label}
\newcommand \bpro {\begin{proof}}
\newcommand \epro {\end{proof}}
\newcommand \bdf {\begin{definition}}
\newcommand \edf {\end{definition}}
\newcommand \bex {\begin{example}}
\newcommand \eex {\end{example}}
\newcommand \bcl {\begin{Claim}}
\newcommand \ecl {\end{Claim}}
\newcommand \brm {\begin{Remark}}
\newcommand \erm {\end{Remark}}

\let\oldmarginpar\marginpar
\renewcommand\marginpar[1]{\-\oldmarginpar[\raggedleft\footnotesize #1]%
{\raggedright\footnotesize #1}}

\begin{document}
 
\begin{frontmatter}

\title{A central-upwind geometry-preserving method for 
\\
hyperbolic conservation laws on the sphere}

\author{Abdelaziz Beljadid$^1$ and Philippe G. LeFloch$^2$}
\address{$^1$ Department of Civil Engineering, University of Ottawa, 161 Louis Pasteur, Ottawa, 
\\
Ontario, K1N6N5, Canada, Email: abelj016@uottawa.ca
\\
$^2$ Laboratoire Jacques-Louis Lions \& Centre National de la Recherche Scientifique, 
\\
Universit\'e Pierre et Marie Curie (Paris 6), 4 Place Jussieu, 75258 Paris, France. 
\\
Email: contact@philippelefloch.org}

\begin{abstract} We introduce a second-order, central-upwind finite volume method for the discretization of nonlinear hyperbolic conservation laws posed on the two-dimensional sphere. The semi-discrete version of the proposed method is based on a technique of local propagation speeds and it is free of any Riemann solver.  The main advantages of our scheme are the high resolution of discontinuous solutions, its low numerical dissipation, and its simplicity for the implementation. The proposed scheme does not use any splitting approach, which is applied in some cases to upwind schemes in order to simplify the resolution of Riemann problems. The semi-discrete form of the scheme is strongly linked to the analytical properties of the nonlinear conservation law and to the geometry of the sphere. The curved geometry is treated here in an analytical way so that the semi-discrete form 
of the proposed scheme is consistent with a geometric compatibility property. Furthermore, the time evolution is carried out by using a total-variation-diminishing Runge-Kutta method. 
A rich family of (discontinuous) stationary solutions is available for the problem under consideration when the flux is nonlinear and foliated (as identified by the author in an earlier work). We present here a series of numerical examples, obtained by considering non-trivial steady state solutions and this leads us to a good validation of the accuracy and efficiency of the proposed central-upwind finite volume method. Our numerical tests confirm the stability of the proposed scheme and clearly show its ability to capture accurately discontinuous steady state solutions to nonlinear
hyperbolic conservation laws posed on the sphere. 
\end{abstract}

\end{frontmatter}



\section{Introduction}
\label{sec1}

Solutions to hyperbolic partial differential equations generally develop discontinuities in finite time, even from smooth initial conditions. Various classes of so-called shock-capturing schemes have been proposed. In particular, 
upwind and central schemes have been used to numerically solve these equations. Generally, it can be stated that the difference between these schemes is that upwind methods use characteristic-related information, while central methods do not. The use of characteristic information in upwind schemes can improve the results but renders these schemes, in some cases, computationally expensive. Central schemes are widely used (see e.g.\cite{Russo2002}) after the pioneering work of Nessyahu and Tadmor \cite{NessyahuTadmor}, where a second-order finite volume central method on a staggered grid in space-time was first proposed.
This strategy leads to high-resolution and the simplicity of the Riemann-solver free method. As observed in Kurganov and Tadmor \cite{KurganovTadmor2000a}, this scheme suffers from excessive numerical viscosity when a small time step is considered.

In order to improve the performance of central schemes, some characteristic information can still be used.
Kurganov et al. \cite{KurganovNP} proposed the central-upwind schemes which are based on information obtained from the local speeds of wave propagation. The central-upwind schemes can be considered as a generalization of central schemes originally developed by Kurganov and Tadmor \cite{KurganovTadmor2000a, KurganovTadmor2000b},  Kurganov and Levy \cite{KurganovLevy}, and Kurganov and Petrova \cite{KurganovPetrova2001}. The central-upwind schemes are simple, since they use no Riemann solvers, and they have proven their effectiveness in multiple studies, as shown in \cite{KurganovPetrova2006, KurganovPetrova2007, KurganovPetrova2008, KurganovPetrova2009, Kurganovetal2007,KurganovLevy2002}. Kurganov and Petrova \cite{KurganovPetrova2005} extended the central-upwind schemes to triangular grids for solving two-dimensional Cartesian systems of conservation laws. Next, Beljadid et al. \cite{BMK} proposed a two-dimensional well-balanced and positivity preserving cell-vertex central-upwind scheme for the computation of shallow water equations with source terms due to bottom topography. 

Several studies have been recently developed for hyperbolic conservation laws posed on curved manifolds. The solutions of conservation laws including the systems on manifolds and on spacetimes were studied in \cite{Rossmanith-2004, MortonSonar2007} and by LeFloch and co-authors \cite{ABL,AmorimLeFlochOkutmustur,BL,BFL} and \cite{PLF}--
\cite{LeFlochOkutmustur}.  More recently, hyperbolic conservation laws for evolving surface were investigated by Dziuk, Kro\" oner and M\"uller \cite{DKM}, Giesselman \cite{Giesselmann}, and Dziuk and Elliott \cite{DE}.
Earlier on, for such problems, Ben-Artzi and LeFloch \cite{BL} and LeFloch and Okutmustur \cite{LeFlochOkutmustur}
established a general well-posedness theory for conservation laws on manifolds. 

Burgers equation provides a simple, yet challenging equation,  which admits discontinuous solutions and it provides a simplified setup for the design and validation of shock-capturing numerical methods. Burgers equation and its generalizations to a curved manifold have been widely used in the physical and mathematical literature. 
In \cite{BLM}, we have used a class of Burgers-type equations on the sphere and adopted the methodology first proposed by Ben-Artzi, Falcovitz, and LeFloch \cite{BFL}, which uses second--order approximations based on generalized Riemann problems. In \cite{BLM}, a scheme was proposed which uses piecewise linear reconstructions based on solution values at the center
 of the computational cells and on values of Riemann solutions at the cell interfaces. A second-order approximation based on a generalized Riemann solver was then proposed, together with a total variation diminishing Runge-Kutta method (TVDRK3) with operator splitting for the temporal integration. 

The finite volume methods developed in \cite{BLM} and \cite{BFL} are strongly linked to the structure of governing equation. The geometric dimensions are considered in an analytical way which leads to discrete forms of schemes that respect exactly the geometric compatibility property. The splitting approach which is used in these schemes simplifies the resolution of the Riemann problem but it increases the computational cost.

In the present study, we propose a new finite volume method which is less expensive in terms of computational cost. This scheme is free of any Riemann solver and does not use any splitting approach, while such a splitting 
is widely used in upwind schemes when one needs to simplify the resolution of Riemann problems. The proposed paper provides the first study of {\sl central-upwind schemes for conservation laws on a curved geometry.} 

Burgers equation and its generalizations will be used in the present paper in order to develop and validate the new finite volume method. We design in full detail a geometry-compatible central-upwind scheme for scalar nonlinear hyperbolic conservation laws on the sphere. This system has a simple appearance but it generates solutions that have a very rich wave structure (due to the curved geometry) and its solutions provide an effective framework for assessing numerical methods. Our goal is to develop and validate a finite volume method which is free of any Riemann problem and is consistent with the geometric compatibility (or divergence free) condition, at the discrete level. As we prove it, the proposed scheme is efficient and accurate for discontinuous solutions and implies only negligible geometric distorsions on the solutions.

An outline of the paper is as follows. In Section \ref{sec2}, the governing equations related to this study are presented.  Section \ref{sec3} is devoted to the derivation of the semi-discrete version of our scheme. In Section \ref{sec4}, the coordinate system and the non-oscillatory reconstruction are described. In Section \ref{sec5}, we present the geometry-compatible flux vectors and some particular steady state solutions as well as confined solutions, which will be used to validate the performance of the proposed method. In Section \ref{sec6}, we demonstrate the high-resolution of the proposed central-upwind scheme thanks to a series of numerical experiments. Finally, some concluding remarks are provided.


\section{Governing equations}
\label{sec2}

We consider nonlinear hyperbolic equations posed on the sphere $\mathbb{S}^{2}$ and based on the flux-vector $F:=F(x, \ubar)$, depending on the function $u(t,x)$ and the space variable $x$. This flux is assumed to satisfy the following geometric compatibility condition: for any arbitrary constant value $\bar{u} \in \RR$, 
\bel{Eq2.1}
 \nabla \cdot (F(\cdot,\bar{u}) ) = 0,
\ee
and that the flux takes the form
\bel{Eq2.2}
F(x,u)=n(x)\wedge \Phi (x,u), 
\ee
where $n(x)$ is the unit normal vector to the sphere and the function $\Phi (x,u)$ is a vector field in $\mathbb{R}^3$, restricted to $\mathbb{S}^{2}$ and defined by
\bel{Eq2.3}
\Phi (x,u)=\nabla h(x,u). 
\ee
Here, $h\equiv h(x,u)$ is a smooth function depending on the space variable $x$ and the state variable $u(t,x)$. Observe that (for instance by Claim 2.2 in \cite{BFL}) the conditions $(\ref{Eq2.2})$ and $(\ref{Eq2.3})$ for the flux vector are sufficient to ensure the validity of the geometric compatibility condition $(\ref{Eq2.1})$.

We are going here to develop and validate a new geometry-preserving central-upwind scheme which approximates solutions to the hyperbolic conservation law 
\bel{Eq2.4}
\del_t u + \nabla \cdot F(x,u)  = 0,   \quad  (x,t)\in \mathbb{S}^{2}\times \mathbb{R}_{+},
\ee
where $\nabla \cdot F$ is the divergence of the vector field $F$. Given any data $u_{0}$ prescribed on the sphere, we consider the following initial condition for the unknown function $u:= u(t,x)$ 
\bel{Eq2.5}
u(x,0)=u_{0}(x), \qquad x \in \mathbb{S}^{2}.
\ee
Equation $\eqref{Eq2.4}$ can be rewritten, using general local coordinates and the index of summation $j$, in the form
\bel{Eq2.6}
\partial _{t}u+\frac{1}{\sqrt{\left | g \right |}}\partial_{j} (\sqrt{\left | g \right |}F^{j}(x,u))=0,
\ee
or
\bel{Eq2.7}
\partial _{t}(\sqrt{\left | g \right |}u)+\partial_{j} (\sqrt{\left | g \right |}F^{j}(x,u))=0,
\ee
where in local coordinates $x=(x^{j})$, the derivatives are denoted by $\partial _{j}:=\frac{\partial }{\partial x^{j}}$, $F^{j}$ are the components of the flux vector and $g$ is the metric.

The conservation law $(\ref{Eq2.4})$ becomes as follows
 \bel{Eq2.8}
\partial _{t} v +\partial_{j} (\sqrt{\left | g \right |}F^{j}(x,v/\sqrt{\left | g \right |}))=0,
\ee
where $v=u\sqrt{\left | g \right |}$. This form will be used in the derivation of the semi-discrete form of the proposed scheme.
For latitude-longitude grid on the sphere, the divergence operator of the flux vector is  
\bel{eq:4.3}
 \nabla \cdot F =  \frac{1}{\cos \phi } (\frac{\partial  }{\partial \phi}(F_{\phi } \cos \phi) +\frac{\partial F_{\lambda }}{\partial \lambda }   ),
\ee
where $F_{\phi }$ and $F_{\lambda}$ are the flux components in the latitude ($\phi$) and longitude $(\lambda)$ directions on the sphere respectively.


\section{Derivation of the proposed method}
\label{sec3}

\subsection{Discretization of the divergence operator}

The derivation of the new central-upwind scheme will be described in detail for the three steps: reconstruction, evolution, and projection. We will develop and give a semi-discrete form of the proposed method for general computational grid 
used to discretize the sphere. We assume the discretization of the sphere $\mathbb{S}^{2}: =\bigcup_{j=1}^{j=N} C_{j}$, where $C_{j}$ are the computational cells with area $\left | C_{j} \right |$. 
We denote by $m_{j}$ the number of cell sides of $ C_{j}$ and by $ C_{j1}, C_{j2},... C_{jm_{j}}$ the neighboring computational cells that share with $C_{j}$ the common sides $(\partial C_{j})_{1}, (\partial C_{j})_{2},...(\partial C_{j})_{m_{j}}$ , respectively. The length of each cell-interface $(\partial C_{j})_{k}$ is denoted by $l_{jk}$. The discrete value of the state variable $u(t,x)$ inside the computational cell $C_{j}$ at a point $G_{j} \in C_{j}$ is denoted by $\mathbf{u}_{j}^{n}$ at step $n$. 
The longitude and latitude coordinates of the suitable point $G_{j}$ to use inside each computational cell $C_{j}$ are presented in Section \ref{sec4.2}. These coordinates should be chosen according to the reconstruction of the state variable $u(t,x)$ over the computational cells used on the sphere.
Finally, we use the notations $\Delta t$ and $t_{n} := n\Delta t$ for the time step and the time at step $n$, respectively.
Note that the development of the first order in time is sufficient to have the exact semi-discrete form of the proposed scheme. The resulting ODE can be numerically solved using a higher-order SSP ODE solver as Runge-Kutta of multistep methods. In the numerical experiments, 
the third-order TVD Runge-Kutta method proposed by Shu and Osher \cite{ShuOsher} is used.

In this section, we will present a general form of the discretization of the divergence operator for general computational grid on the sphere. The approximation of the flux divergence can be written using the divergence theorem as

\bel{Eq3.1}
\begin{aligned}
&[\nabla \cdot F(x,u)]^{approx}=       \frac{I_{j}}{\left | C_{j} \right |},   \\
&I_{j}=[\oint _{\partial C_{j}}F(x,u)\cdot\nu(x)ds]^{approx},
\end{aligned}
\ee
where $\nu(x)$ is the unit normal vector to the boundary $\partial C_{j}$ of the computational cell $C_{j}$ and $ds$ is the infinitesimal length along $\partial C_{j}$. 

The scalar potential function $h$ is used to obtain the following approximation along each side of the computational cell $C_{j}$. 

\begin{Claim}
\label{Claim3.1}
For a three-dimensional flux $\Phi(x,u)$ given by $\eqref{Eq2.3}$, where $h\equiv h(x,u)$ is a smooth function in the neighborhood of the sphere $\mathbb{S}^{2}$, the total approximate flux through 
the cell interface $e$ is given by
\bel{Eq3.2}
\oint _{e^{1}}^{e^{2}}F(x,u)\cdot\nu(x)ds=-(h(e^{2},u_{j})-h(e^{1},u_{j})),
\ee
where $e^{1}$ and $e^{2}$ are the  initial and final endpoints of the side $e$ using the sense of integration and $u_{j}$ is the estimate value of the variable $u$ along the side $e$.
\end{Claim}


Namely, the flux vector is written in the form 
$F(x,u)=n(x)\wedge \Phi(x,u)$ and we can derive the approximation of the integral along each cell side of $C_{j}$ 
\bel{Eq3.3}
\begin{aligned}
&\oint _{e^{1}}^{e^{2}}F(x,u)\cdot\nu(x) d s=\oint _{e^{1}}^{e^{2}}(n(x)\wedge \Phi(x,u))\cdot\nu(x)d s \\
&= -\oint _{e^{1}}^{e^{2}}\Phi(x,u)\cdot(n(x)\wedge \nu(x))d s=
-\oint _{e^{1}}^{e^{2}} \nabla h(x,u)\cdot \tau (x)d s\\
&=-\oint _{e^{1}}^{e^{2}} \nabla _{\partial C_{j}}h(x,u)d s=-(h(e^{2},u_{j})-h(e^{1},u_{j})),
\end{aligned}
\ee
where $\tau(x)$ is the unit vector tangent to the boundary $\partial C_{j}$.



\begin{remark}
Using the discrete approximations based on Claim $\ref{Claim3.1}$, if a constant value of the state variable $u(t,x)=u_{j}=\bar{u}$ is considered one obtains
\bel{Eq3.4}
\begin{aligned}
&[\nabla \cdot F(x,u)]^{approx}= \frac{1}{\left | C_{j} \right |}[\oint _{\partial C_{j}}F(x,u)\cdot\nu(x)ds]^{approx} \\
&=-\sum_{e\in \partial C_{j}} (h(e^{2},\bar{u})-h(e^{1},\bar{u}))=0.
\end{aligned}
\ee
This confirms that the discrete approximation of the divergence operator respects the divergence free condition which is the geometric requirement that the proposed scheme should satisfy. 
\end{remark}

\subsection{Reconstruction method}
\label{Sec3.2.2}
In the following, we will present the reconstruction of the proposed central-upwind scheme. The semi-discrete form of the proposed scheme for Equation $(\ref{Eq2.4})$ will be derived by using the approximation of the cell averages of the solution. At each time $t=t_{n} $ the computed solution is
\bel{Eq3.5}
\mathbf{u}_{j}^{n}\approx \frac{1}{\left | C_{j} \right |}\int_{C_{j}}^{}u(x,t_{n} ) dV_{g},
\ee
where $dV_{g} = \sqrt g dx^{1}dx^{2}$.

The discrete values $\mathbf{u}_{j}^{n}$ of the solution at time $t=t_{n} $ are used to construct a conservative piecewise polynomial function with possible discontinuities at the interfaces of the computational cells $C_{j}$
\bel{Eq3.6}
\tilde{u}^{n}(x )=\sum_{j}^{}w_{j}^{n}(x)\chi_{j}(x),
\ee
where $w^{n}_{j}(x)$ is a polynomial in two variables ($\lambda$ and $\phi$), and $\chi_{j}$ is the characteristic function which is defined using the Kronecker symbol $\delta _{jk}$ and for any point of spatial coordinate $x$ inside the computational cell $C_{k}$ we consider $\chi_{j}(x):=\delta _{jk}$.

The maximum of the directional local speeds of propagation of the waves at the \textit{kth} interface inward and outward of the computational cell $C_{j}$ are denoted by $a_{jk}^{in}$ and $a_{jk}^{out}$, respectively. 
When the solution evolves over a time step $\Delta t$, the discontinuities move inward and outward the \textit{kth} interface of the computational cell $C_{j}$ with maximum distances $a_{jk}^{in} \Delta t$ and $a_{jk}^{out} \Delta t$, respectively. These distances of propagation are used at the computational cells to delimit different areas in which the solution still smooth and the areas in which the solution may not be smooth when it evolves from the time level $t_{n} $ to $t_{n+1}$. 

\begin{figure}[htbp]
\begin{centering}
\includegraphics[width=11cm,height=8cm,angle=0]{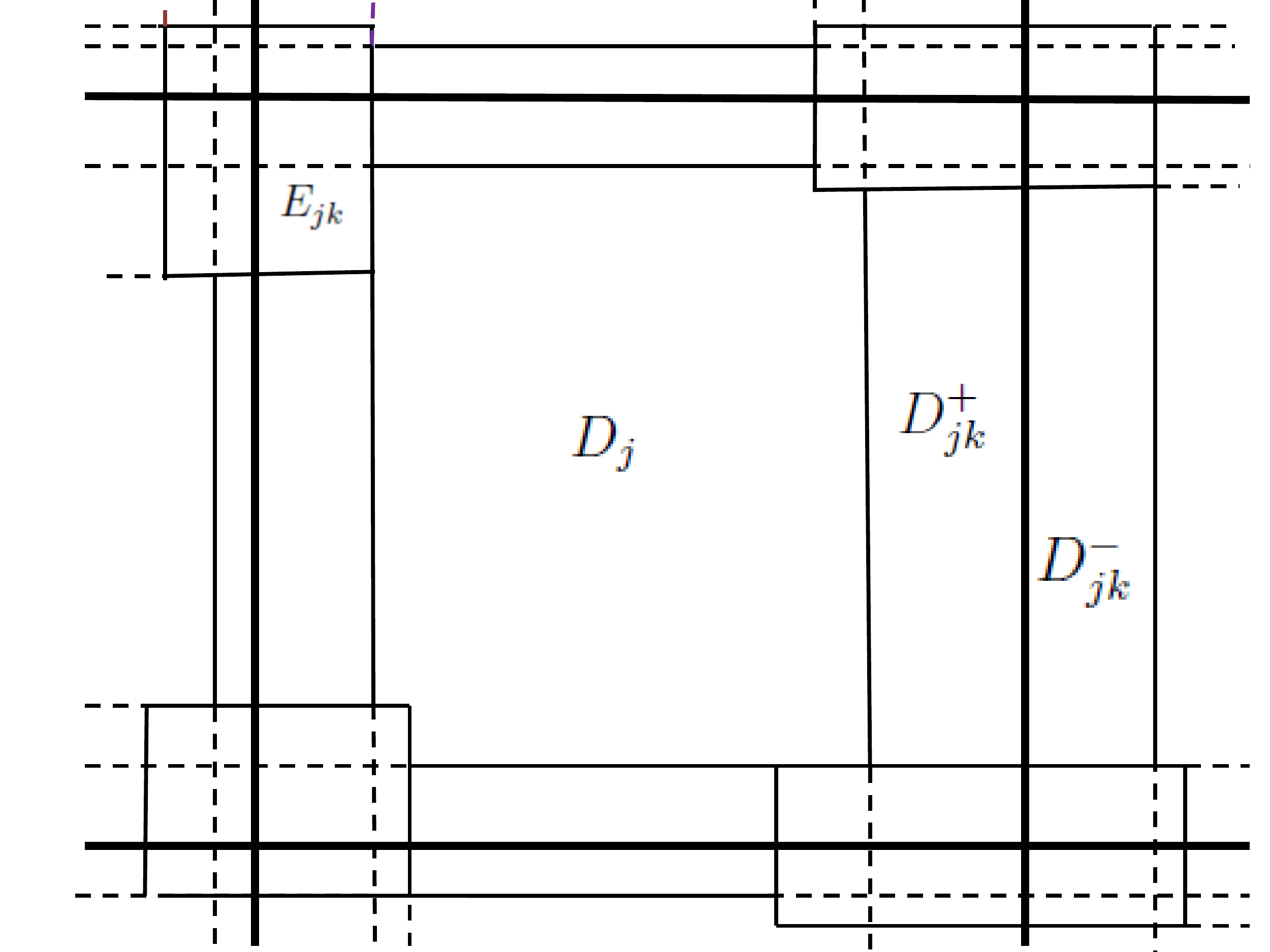} 
\par\end{centering}
\caption{Schematic view of the decomposition of the control volume}
\label{Fig1} 
\end{figure}

We define the domain $D_{j}$ as the part inside the cell $C_{j}$ in which the solution still smooth, see Figure \ref{Fig1}. Two other types of domains are defined, the first type includes the rectangular domains $D_{jk}$, $k=1,2...m_{j}$, along each side of $C_{j}$ of width $ (a_{jk}^{out}+a_{jk}^{in})\Delta t$ and length $l_{jk} +O(\Delta t)$ and the second type includes the domains denoted by $E_{jk}$, $k=1,2...m_{j}$, around the cell vertices of computational cells. These domains are decomposed into two sub-domains $D_{jk}=D^{+}_{jk} \cup  D^{-}_{jk}$ and $E_{jk}=E^{+}_{jk} \cup  E^{-}_{jk}$, where the sub-domains with the superscripts 
plus signs ``$+$'' and minus signs ``$-$'' are the domains inside and outside of the cell $C_{j}$, respectively. For purely geometrical reasons, the areas of the three types of sub-domains are of orders $\left | D_{j} \right |= O(1)$, $\left | D_{jk} \right |= O(\Delta t)$ and $\left | E_{jk} \right |= O(\Delta t^{2})$. 

We consider the projection of the flux vector $\tilde{F}$ according to the orthogonal to the \textit{kth} cell interface $(\partial C_{j})_{k}$
\bel{Eq3.7}
\begin{aligned}
f_{jk}=N_{jk}\cdot \tilde{F},
\end{aligned}
\ee
where $N_{jk}$ is the unit normal vector to the cell interface $(\partial C_{j})_{k}$ and $\tilde{F}$ has the components $\sqrt{ g }F^{j}(x,v/\sqrt{ g})$ which are used in Equation (\ref{Eq2.8}).

The one-sided local speeds of propagation of the waves at the \textit{kth} cell interface $(\partial C_{j})_{k}$,   inward and outward the  computational cell $C_{j}$,  are estimated by
\bel{Eq3.8}
\begin{aligned}
&a^{out}_{jk}=\max \{\frac{\partial f_{jk}}{\partial v}(M_{jk},u_{j}(M_{jk})),\frac{\partial f_{jk}}{\partial v}(M_{jk},u_{jk}(M_{jk})),0 \},\\
&a^{in}_{jk}=-\min \{\frac{\partial f_{jk}}{\partial v}(M_{jk},u_{j}(M_{jk})),\frac{\partial f_{jk}}{\partial v}(M_{jk},u_{jk}(M_{jk})),0 \},\\
\end{aligned}
\ee
where $u_{j}(M_{jk})$ is the value of the state variable $u$ at the midpoint $M_{jk}$ of $(\partial C_{j})_{k}$, which is obtained from the non-oscillatory reconstruction for the computational cell $C_{j}$ and $u_{jk}(M_{jk})$ is the value 
of $u$ at the same point $M_{jk}$ using the non-oscillatory reconstruction for the neighboring cell $C_{jk}$.

\subsection{Evolution and projection steps}

The computed cell averages $\mathbf{\bar{u}}_{j}^{n+1}$ of the numerical solution at time step $t_{n+1}$ over the computational cells $C_{j}$ are used to obtain the piecewise linear reconstruction $\tilde{w}^{n+1}$ which should satisfy the following conservative 
requirement
\begin{equation}
\mathbf{\bar{u}}_{j}^{n+1}= \frac{1}{\left | C_{j} \right |}\int_{C_{j}}^{}\tilde{w}^{n+1}(x)dV_{g}.
\qquad\label{Eq3.9}
\end{equation}

The average of the function $\tilde{w}^{n+1}$ over the domain $D_{j}$ is denoted by $\bar{w}^{n+1}(D_{j})$ 
\begin{equation}
\bar{w}^{n+1}(D_{j})= \frac{1}{\left | D_{j} \right |}\int_{D_{j}}^{}\tilde{w}^{n+1}(x) dV_{g}.
\qquad\label{Eq3.10}
\end{equation}

Note that it is possible to derive the fully discrete form of the proposed scheme but it is impractical to use and for simplicity, we will develop the semi-discrete form of the scheme. The ODE for approximating the cell averages of the solutions is derived by tending the time step $ \Delta t$ to zero. This eliminates some terms because of their orders and we keep the more consistent terms
\begin{equation}
\begin{aligned}
&\frac{d\mathbf{\bar{u}}_{j} }{dt}(t_{n} )=  \underset{\Delta t\rightarrow 0}{\lim}\frac{\mathbf{\bar{u}}_{j}^{n+1}-\mathbf{\bar{u}}_{j}^{n}}{\Delta t} \\
 &=\underset{\Delta t\rightarrow 0}{\lim}  \frac{1}{\Delta t}     [   \frac{1}{|  C_{j} |}   \int_{D_{j}}\tilde{w}^{n+1}(x)dV_{g}+\frac{1}{|  C_{j} |} \sum_{k=1}^{m_{j}} \int_{D^{+}_{jk}}\tilde{w}^{n+1}(x)dV_{g}  \\
&+ \frac{1}{|  C_{j} |} \sum_{k=1}^{m_{j}} \int_{E^{+}_{jk}}\tilde{w}^{n+1}(x)dV_{g}-\mathbf{\bar{u}}_{j}^{n}  ].
 \end{aligned}
\qquad\label{Eq3.11}
\end{equation}

Since the areas of domains $E_{jk}$ with $k=1,2,..., m_{j}$ are of order $\Delta t^{2}$ we obtain 
\begin{equation}
\int_{E^{+}_{jk}}\tilde{w}^{n+1}(x)dV_{g}=O(\Delta t^{2}).
\qquad\label{Eq12}
\end{equation}

This approximation allows us to deduce that the third term on the right-hand side of Equation $(\ref{Eq3.11})$ is of order $\Delta t ^{2}$ and the result for the limit of this term vanishes for the ODE.

The second term in Equation $(\ref{Eq3.11})$, in which we use the rectangular domains $D_{jk}^{+}$, will be estimated by using the assumption that the spatial derivatives of $\tilde{w}^{n+1}$ are bounded independently of $\Delta t$. Under this assumption the following Claim gives an estimation of this term with an error of order $\Delta t^{2} $ for each $k\in [1,m_{j}]$.
\begin{Claim}
\label{Claim3.2}
Consider the reconstruction given by $(\ref{Eq3.6})$, its evolution $\tilde{w}^{n+1}$ over the global domain, and the definitions given in Section \ref{Sec3.2.2} for the domains $D_{jk}$ and $ D^{+}_{jk}$.  If we assume that the spatial derivatives of $\tilde{w}^{n+1}$ are bounded independently of $\Delta t$, then 
\begin{equation}
\int_{D^{+}_{jk}}\tilde{w}^{n+1}(x)dV_{g}= | D_{jk} ^{+}|  \bar{w}^{n+1}(D_{jk})+O(\Delta t^{2}  ).
\qquad\label{Eq3.13}
\end{equation}
\end{Claim}

The proof is as follows. 
It is obvious that for the case $\left | D_{jk}^{+} \right |=0$ or $\left | D_{jk}^{-} \right |=0$ equation $\eqref{Eq3.13}$ is valid. 
We assume that $\left | D_{jk}^{+} \right |    \left | D_{jk}^{-} \right | \neq 0$ and we consider 
$$R=\int_{D^{+}_{jk}}\tilde{w}^{n+1}(x)dV_{g}- | D_{jk} ^{+}|  \bar{w}^{n+1}(D_{jk})$$
We have
\begin{equation}
\begin{aligned}
 &R=\int_{D^{+}_{jk}}\tilde{w}^{n+1}(x)dV_{g}-\frac{\left | D^{+}_{jk} \right |}{\left | D_{jk} \right |} (\int_{D^{+}_{jk}}\tilde{w}^{n+1}(x)dV_{g}+\int_{ D^{-}_{jk} }\tilde{w}^{n+1}(x)dV_{g}) \\
&=\frac{|D^{+}_{jk}|}{| D_{jk}|} [ \frac{|D^{-}_{jk}|}{|D^{+}_{jk}|} \int_{D^{+}_{jk}}\tilde{w}^{n+1}(x)dV_{g} - \int_{D^{-}_{jk}}\tilde{w}^{n+1}(x)dV_{g} ] \\
&=\frac{ | D^{+}_{jk}  |}{ |D_{jk} |} [\frac{a^{out}_{jk}}{ a^{in}_{jk}}\int_{-a^{in}_{jk}\Delta t}^{0}  \tilde{w}^{n+1}(s)\tilde{l}_{jk}ds - \int_{0}^{a^{out}_{jk}\Delta t}\tilde{w}^{n+1}(s)\tilde{l}_{jk}ds   ],
 \end{aligned}
\qquad\label{Eq15}
\end{equation}
where $\tilde{l}_{jk}$ is the length of the domain $D_{jk}$ and $s$ is a variable according to the orthogonal outward axis to the \textit{kth} cell interface, see Figures \ref{Fig1} and \ref{Fig3}.    

One obtains after the change of variable in the first integral of the last equality in $(\ref{Eq15})$
$$R=\frac{| D^{+}_{jk} |}{|D_{jk}|}\tilde{ l}_{jk}  \int_{0}^{a^{out}_{jk}\Delta t} (\tilde{w}^{n+1}(-\frac{a^{in}_{jk}}{a^{out}_{jk}}s)-\tilde{w}^{n+1}(s)   )ds.$$
Using the mean value theorem to the function $\tilde{w}^{n+1}$ we obtain
$$R=-\frac{| D^{+}_{jk} |}{|D_{jk}|}\tilde{ l}_{jk}\int_{0}^{a^{out}_{jk}\Delta t} \frac{a^{in}+a^{out}}{a^{out}}s\frac{\partial \tilde{w}^{n+1}}{\partial s}(c_{s})ds,$$
 where $c_{s}\in [min(s,-s a^{in}_{jk}/a^{out}_{jk}),max(s,-s a^{in}_{jk}/a^{out}_{jk})]$. 

We denote by $M$ the upper bound value of the spatial derivative of the function $\tilde{w}^{n+1}$ over the domain $D_{jk}$. Therefore we obtain
$$|R|\leq M l \frac{| D^{+}_{jk} |}{|D^{-}_{jk}|}\int_{0}^{a^{out}_{jk}\Delta t}sds= \frac{ M l}{2} |  D^{+}_{jk} ||  D^{-}_{jk} |.$$
Since $\tilde{l}_{jk}=l_{jk}+O(\Delta t) $, and both the areas $|  D^{+}_{jk}|$ and  $|  D^{-}_{jk} |$ are of order $\Delta t$ we obtain $R=O(\Delta t^{2})$.
This completes the proof. 

\


Using Equation $(\ref{Eq3.13})$ in Claim $\ref{Claim3.2}$ one obtains
\begin{equation}
\begin{aligned}
&  \frac{1}{|  C_{j} |} \sum_{k=1}^{m_{j}} \int_{D^{+}_{jk}}\tilde{w}^{n+1}(x)dV_{g}=\frac{1}{| C_{j}|}  \sum_{k=1}^{m_{j}} |D_{jk}^{+} |\bar{w}^{n+1} (D_{jk})+O(\Delta t^{2})  \\
&= \frac{\Delta t}{| C_{j}|} \sum_{k=1}^{m_{j}}a_{jk}^{in}(l_{jk}+O(\Delta t))\bar{w}^{n+1}(D_{jk})+O(\Delta t^{2}).
\end{aligned}
\qquad\label{Eq3.14}
\end{equation}

Therefore Equation $(\ref{Eq3.11})$ can be written as
\begin{equation}
\frac{d\mathbf{\bar{u}}_{j} }{dt}(t_{n} )= \underset{\Delta t\rightarrow 0}{\lim} \frac{1}{\Delta t} [ \frac{|D_{j}|}{|C_{j}|  }\bar{w}^{n+1}(D_{j})-\mathbf{\bar{u}}_{j}^{n}]+\sum_{k=1}^{m_{j}}  \underset{\Delta t\rightarrow 0}{\lim} \frac{|D_{jk}^{+}|}{\Delta t |C_{j}| } \bar{w}^{n+1}(D_{jk}),
\qquad\label{Eq3.15}
\end{equation}
where 
\begin{equation}
\bar{w}^{n+1}(D_{jk})=\frac{1}{\left | D_{jk} \right |}\int_{D_{jk}} \tilde{w}^{n+1}  (x) dV_{g}.
\qquad\label{Eq3.16}
\end{equation}

In order to derive the semi-discrete form of the proposed scheme from Equation $(\ref{Eq3.15})$, one needs to compute the average values $\bar{w}^{n+1}(D_{jk})$ and $\bar{w}^{n+1}(D_{j})$. 
To compute $\bar{w}^{n+1}(D_{jk})$, Equation $(\ref{Eq2.4})$ is integrated over the space-time  control volume $D_{jk} \times [t_{n},t_{n+1}]$. 
After integration by parts and applying the divergence theorem to transform the surface integral of the divergence operator to the boundary integral and using the approximation (\ref{Eq3.2}) of the flux through the cell interfaces, the following equations are obtained
\begin{equation}
\begin{aligned}
& \bar{w}^{n+1}(D_{jk})=\frac{1}{|D_{jk}|}[\int_{D_{jk}^{+}}w_{j}^{n}(x)dV_{g}+\int_{D_{jk}^{-}}w_{jk}^{n}(x)dV_{g}]
  \\
& -\frac{1}{|D_{jk}|}\int_{t_{n} }^{t_{n+1}}\int_{D_{jk}} \nabla \cdot F(x,u)dV_{g}            ,
 \end{aligned}
\qquad\label{Eq3.17}
\end{equation}
and
\begin{equation}
\begin{aligned}
& \int_{D_{jk}}\nabla \cdot F(x,u)dV_{g}= [\int_{\partial D_{jk}} F(x,u)\cdot \nu(x) ds]^{approx}= \sum_{i=1}^{i=4}\int_{(\partial D_{jk})_{i}}F(x,u) \cdot \nu(x) ds \\
&=-[-h(e^{2}_{jk},u_{j}(M_{jk}))+h(e^{1}_{jk},u_{j}(M_{jk})) +  h(e^{2}_{jk},u_{jk}(M_{jk}))-h(e^{1}_{jk},u_{jk}(M_{jk}))     ]+O(\Delta t),
 \end{aligned}
\qquad\label{Eq3.18}
\end{equation}
where $(\partial D_{jk})_{i}$, $i=1,2,3,4$, are the four edges of the rectangular domain $D_{jk}$, $e^{2}_{jk}$ and $e^{1}_{jk}$ are the initial and final endpoints 
of the cell interface $(\partial C_{j})_{k}$, and as mentioned before $w_{j}^{n}$ and $w_{jk}^{n}$ are the piecewise polynomial reconstructions in the computational cells $C_{j}$ and $C_{jk}$ respectively at time $t_{n}$.

The term on the right-hand side of Equation $(\ref{Eq3.18})$ of order $O(\Delta t)$ corresponds to the global result of the integration along the two edges of the domain $D_{jk}$ having the length $(a^{int}_{jk}+a^{out}_{jk}) \Delta t $, and 
the rest of the integration due to the difference between the length of the domain $D_{jk}$ and the length of the cell interface $(\partial C_{j})_{k}$. 

In order to compute the spatial integrals in Equation $\eqref{Eq3.17}$, the Gaussian quadrature can be applied. In our case, the midpoint rule is used for simplicity
\begin{equation}
\int_{D^{+}_{jk}}w^{n}_{jk}dV_{g}+\int_{D^{-}_{jk}}w^{n}_{jk}dV_{g}\approx  l_{jk} 
\Delta t    [a^{in}_{jk}u_{j}(M_{jk})+a^{out}_{jk}u_{jk}(M_{jk})].
\qquad\label{Eq3.19}
\end{equation}
Equations \eqref{Eq3.17}, \eqref{Eq3.18} and \eqref{Eq3.19} lead to  

\begin{equation}
\begin{aligned}
& \lim_{\Delta t\rightarrow 0} \bar{w}^{n+1}(D_{jk})=\frac{l_{jk}}{a^{in}_{jk}+a^{out}_{jk}}
[a^{in}_{jk}u_{j}(M_{jk})+a^{out}_{jk}u_{jk}(M_{jk})]\\
& +\frac{1}{a^{in}_{jk}+a^{out}_{jk}}[-h(e^{2}_{jk},u_{j}(M_{jk}))+h(e^{1}_{jk},u_{j}(M_{jk})) +  h(e^{2}_{jk},u_{jk}(M_{jk}))-h(e^{1}_{jk},u_{jk}(M_{jk}))  ].
\qquad\label{Eq3.20}
\end{aligned}
\end{equation}
Therefore
\begin{equation}
\begin{aligned}
& \lim_{\Delta t\rightarrow 0}\sum_{k=1}^{m_{j}}\frac{\left | D_{jk}^{+} \right |}{\Delta t \left | C_{j} \right |}\bar{w}^{n+1}(D_{jk})=\sum_{k=1}^{m_{j}}\frac{a^{in}_{jk}l_{jk}}{\left | C_{j} \right |(a^{in}_{jk}+a^{out}_{jk})}
[a^{in}_{jk}u_{j}(M_{jk})+a^{out}_{jk}u_{jk}(M_{jk})] \\
&+\sum_{k=1}^{m_{j}}\frac{a^{in}_{jk}}{\left | C_{j} \right | (a^{in}_{jk}+a^{out}_{jk})}[-h(e^{2}_{jk},u_{j}(M_{jk}))+h(e^{1}_{jk},u_{j}(M_{jk})) +  h(e^{2}_{jk},u_{jk}(M_{jk}))-h(e^{1}_{jk},u_{jk}(M_{jk}))  ].
\qquad\label{Eq3.21}
\end{aligned}
\end{equation}
Now, the average value $\bar{w}^{n+1}(D_{j})$ will be computed. Equation $\eqref{Eq2.4}$ is integrated over the space-time control volume $D_{j} \times [t_{n},t_{n+1}]$ and after integration by parts, using the divergence theorem to transform the surface integral to boundary integral and Equation (\ref{Eq3.2}) one obtains
\begin{equation}
\begin{aligned}
& \bar{w}^{n+1}(D_{j})=\frac{1}{|D_{j}|}\int_{D_{j}}w_{j}^{n}dV_{g}-\frac{1}{|D_{j}|}\int_{t_{n}}^{t_{n+1}}\int_{D_{j}}\nabla \cdot F(x,u)dV_{g} \\
&=\frac{1}{|D_{j}|}\int_{D_{j}}w_{j}^{n}dV_{g}-\frac{\Delta t}{\left | D_{j} \right |}
\sum_{k=1}^{m_{j}}[-h(e^{2}_{jk},u_{j}(M_{jk})) + h(e^{1}_{jk},u_{j}(M_{jk}))].
 \end{aligned}
\qquad\label{Eq3.22}
\end{equation}
Using the previous equality we obtain 

\begin{equation}
\begin{aligned}
& \frac{1}{\Delta t}[\frac{\left | D_{j} \right |}{\left | C_{j} \right |} \bar{w}^{n+1}(D_{j})-\bar{u}^{n}_{j} ]=\\
&\frac{1}{\Delta t}\{\frac{1}{|C_{j}|}\int_{D_{j}}w_{j}^{n}dV_{g}-\frac{\Delta t}{\left | C_{j} \right | }\sum_{k=1}^{m_{j}}[-h(e^{2}_{jk},u_{j}(M_{jk})) + h(e^{1}_{jk},u_{j}(M_{jk}))]- \bar{u}^{n}_{j}\},
 \end{aligned}
\qquad\label{Eq3.23}
\end{equation}
which leads to

\begin{equation}
\begin{aligned}
& \lim_{\Delta t\rightarrow 0}\frac{1}{\Delta t}[\frac{\left | D_{j} \right |}{\left | C_{j} \right |} \bar{w}^{n+1}(D_{j})-\bar{u}^{n}_{j} ]=\\
& -\frac{1}{\left | C_{j} \right |}\sum_{k=1}^{m_{j}}a^{in}_{jk}l_{jk}u_{j}(M_{jk})-\frac{1}{\left | C_{j} \right | }\sum_{k=1}^{m_{j}}[-h(e^{2}_{jk},u_{j}(M_{jk})) + h(e^{1}_{jk},u_{j}(M_{jk}))].
 \end{aligned}
\qquad\label{Eq3.24}
\end{equation}
Equations $\eqref{Eq3.21}$ and $\eqref{Eq3.24}$ are used together to obtain the following semi-discrete form 

\begin{equation}
\begin{aligned}
& \frac{d\mathbf{\bar{u}}_{j} }{dt}=-\frac{1}{\left | C_{j} \right |}\sum_{k=1}^{m_{j}}a^{in}_{jk}l_{jk}u_{j}(M_{jk})-\frac{1}{\left | C_{j} \right | }\sum_{k=1}^{m_{j}}[-h(e^{2}_{jk},u_{j}(M_{jk})) + h(e^{1}_{jk},u_{j}(M_{jk}))]+ \\
&\sum_{k=1}^{m_{j}}\frac{a^{in}_{jk}l_{jk}}{\left | C_{j}\right | (a^{in}_{jk}+a^{out}_{jk})}[a^{in}_{jk}u_{j}(M_{jk})+ a^{out}_{jk}u_{jk}(M_{jk}) ]+ \\
& \sum_{k=1}^{m_{j}}\frac{a^{in}_{jk}}{\left | C_{j}\right | (a^{in}_{jk}+a^{out}_{jk})}[-h(e^{2}_{jk},u_{j}(M_{jk})) + h(e^{1}_{jk},u_{j}(M_{jk}))+ h(e^{2}_{jk},u_{jk}(M_{jk})) - h(e^{1}_{jk},u_{jk}(M_{jk})).
 \end{aligned}
\qquad\label{Eq3.25}
\end{equation}

This equation can be rewritten in the following form
\begin{equation}
\begin{aligned}
&\frac{d\mathbf{\bar{u}}_{j}}{dt}=\frac{1}{\left | C_{j} \right |}\sum_{k=1}^{m_{j}}\frac{a^{in}_{jk}a^{out}_{jk}l_{jk}}{a^{in}_{jk}+a^{out}_{jk}}(u_{jk}(M_{jk})-u_{j}(M_{jk}))
+\frac{a^{in}_{jk}a^{out}_{jk}}{\left | C_{j} \right | (a^{in}_{jk}+a^{out}_{jk}) }\{a^{in}_{jk}[h(e^{2}_{jk},u_{j}(M_{jk}))\\
&-h(e^{1}_{jk},u_{j}(M_{jk}))]+
a^{out}_{jk}[h(e^{2}_{jk},u_{jk}(M_{jk}))-h(e^{1}_{jk},u_{jk}(M_{jk}))]\},
 \end{aligned}
\qquad\label{Eq3.26}
\end{equation}
which can be rewritten as follows

\begin{equation}
\begin{aligned}
\frac{d\mathbf{\bar{u}}_{j}}{dt}=-\frac{1}{\left | C_{j} \right |}\sum_{k=1}^{m_{j}}\frac{a^{in}_{jk}\mathbf{H}(u_{jk}(M_{jk}))+a^{out}_{jk}\mathbf{H}(u_{j}(M_{jk}))}{a^{in}_{jk}+a^{out}_{jk}}+\frac{1}{\left | C_{j} \right |}\sum_{k=1}^{m_{j}}\frac{a^{in}_{jk}a^{out}_{jk}l_{jk}}{a^{in}_{jk}+a^{out}_{jk}}[u_{jk}(M_{jk})-u_{j}(M_{jk})],
 \end{aligned}
\qquad\label{Eq3.27}
\end{equation}
where $\mathbf{H}(u_{j}(M_{jk}))$ and $\mathbf{H}(u_{jk}(M_{jk}))$ are given by

\begin{equation}
\begin{aligned}
& \mathbf{H}(u_{j}(M_{jk}))=-[h(e_{jk}^{2},u_{j}(M_{jk}))-h(e_{jk}^{1},u_{j}(M_{jk}))] \\
&\mathbf{H}(u_{jk}(M_{jk}))=-[h(e_{jk}^{2},u_{jk}(M_{jk}))-h(e_{jk}^{1},u_{jk}(M_{jk}))].
 \end{aligned}
\qquad\label{Eq3.28}
\end{equation}

The function $\mathbf{H}$ is defined in the form (\ref{Eq3.28}) in order to be consistent with the total approximate flux through the cell interface 
as presented by Equation (\ref{Eq3.2}) in Claim \ref{Claim3.1}.\\

\smallskip
\noindent
{\bf Remark 2.} If the value of $a^{\rm in}_{jk}+a^{\rm out}_{jk}$ in Equation \eqref{Eq3.27} is zero or very close to zero (smaller than
$10^{-8}$ in our numerical experiments), we avoid division by zero or by a very small number using the following approximations
\begin{equation}
\begin{aligned}
&\frac{a^{in}_{jk}\mathbf{H}(u_{jk}(M_{jk}))+a^{out}_{jk}\mathbf{H}(u_{j}(M_{jk}))}{a^{in}_{jk}+a^{out}_{jk}}\approx \frac{ 1}{2} [\sum_{k=1}^{m_{j}}\mathbf{H}(u_{j}(M_{jk}))+ \sum_{k=1}^{m_{j}}\mathbf{H}(u_{jk}(M_{jk}))],\\
&\frac{a^{in}_{jk} a^{out}_{jk}}{\left | C_{j} \right |(a^{in}_{jk}+a^{out}_{jk})}\sum_{k=1}^{m_{j}}l_{jk}[u_{jk}(M_{jk})-u_{j}(M_{jk})]\approx 0.
\end{aligned}
\qquad\label{Eq3.29}
\end{equation}

These approximations are obtained using similar extreme distances of the propagation of the waves at the cell interface inward and outward the computational cell to define the domains $D_{j}$, $D_{jk}$ and $E_{jk}$.  

The semi-discretization \eqref{Eq3.27} and \eqref{Eq3.28} is a system of ODEs, which has to be integrated in time using an accurate and stable temporal scheme. In our numerical examples reported in Section \ref{sec6}, we used the third-order 
total variation diminishing Runge-Kutta method.

\subsection{The geometry-compatible condition}

In the semi-discrete form \eqref{Eq3.27} and \eqref{Eq3.28} of the proposed scheme, if we consider a constant value of the function $u\equiv\bar{u}$, the second term in the right-hand side of Equation \eqref{Eq3.27} vanishes. For this constant function we obtain for each interface cell $k$  
\begin{equation}
\begin{aligned}
u_{j}(M_{jk})=u_{jk}(M_{jk})=\bar{u},
\end{aligned}
\qquad\label{Eq3.31}
\end{equation}
and
\begin{equation}
\begin{aligned}
\mathbf{H}(u_{j}(M_{jk}))=\mathbf{H}(u_{jk}(M_{jk}))
\end{aligned}
\qquad\label{Eq3.32}
\end{equation}

The first term in the right-hand side of Equation \eqref{Eq3.27} becomes

\begin{equation}
\begin{aligned}
-\frac{1}{\left | C_{j} \right |}\sum_{k=1}^{m_{j}}\frac{a^{in}_{jk}\mathbf{H}(u_{jk}(M_{jk}))+a^{out}_{jk}\mathbf{H}(u_{j}(M_{jk}))}{a^{in}_{jk}+a^{out}_{jk}}=-\frac{1}{\left | C_{j} \right |}\sum_{k=1}^{m_{j}}\mathbf{H}(u_{j}(M_{jk}))
\end{aligned}
\qquad\label{Eq3.33}
\end{equation}

Since we have

\begin{equation}
\begin{aligned}
\sum_{k=1}^{m_{j}}\mathbf{H}(u_{j}(M_{jk}))=\sum_{k=1}^{m_{j}}\mathbf{H}(u_{jk}(M_{jk}))=-\sum_{k=1}^{m_{j}}[h(e^{2}_{jk},\bar{u})-h(e^{1}_{jk},\bar{u})   ]=0, 
\end{aligned}
\qquad\label{Eq3.34}
\end{equation}
we conclude that the first term on the right-hand side of Equation \eqref{Eq3.27} will be canceled which confirms that the proposed scheme respects the geometry-compatibility condition.\\


\begin{remark}
In the proposed central-upwind finite volume method, the midpoint rule was used to compute the spatial integrals. In order to improve the accuracy of the proposed 
scheme, the Gaussian quadrature can be used instead of the midpoint rule. The Gaussian quadrature will not have any impact on the geometry-compatibility condition of the proposed scheme.
\end{remark} 


\section{Formulation using the latitude-longitude grid on the sphere} 
\label{sec4}

\subsection{Computational grid on the sphere}

The geometry-compatible scheme was developed in the previous section for scalar nonlinear hyperbolic conservation laws using general grid on the sphere.  However, in order 
to prevent oscillations an appropriate piecewise linear reconstruction should be proposed according to the computational grid used in the proposed method. 
In the following, we will describe the computational grid and the non-oscillatory piecewise linear reconstruction used in our numerical experiments.

\begin{figure}[h]
\begin{centering}
\includegraphics[width=12cm,height=9.6cm,angle=0]{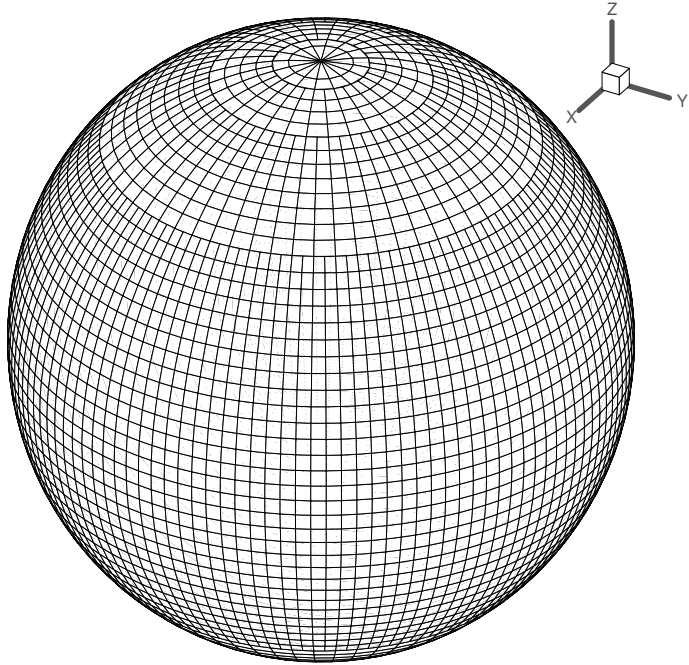} 
\par\end{centering}
\caption{Type of grid used on the sphere}
\label{Fig2}
\end{figure}

The position of each point on the sphere can be represented by its longitude $\lambda \in [0, 2 \pi ]$ and its latitude $\phi \in [-\pi/2, \pi/2 ]$. 
The grid considered in our numerical examples is shown in Figure \ref{Fig2}. The coordinates are singular at the south and north poles, corresponding to $\phi=-\pi/2$ and $\phi= \pi/2$, respectively.
The Cartesian coordinates are denoted by $x=(x_{1},x_{2},x_{3})^{T} \in \mathbb{R}^{3} $ for a standard orthonormal basis vectors $\mathbf{i_{1}}, \mathbf{i_{2}}$, and $\mathbf{i_{3}}$. 

The unit tangent vectors in the directions of longitude and latitude at each point $x$ on the sphere of coordinates $(\lambda,\phi )$ are given as follows
\bel{eq:4.1}
\begin{aligned}
&i_{\lambda }=- \sin\lambda \mathbf{i_{1}}+ \cos\lambda\mathbf{i_{2}},                              \\
& i_{\phi }=- \sin\phi \cos\lambda  \mathbf{i_{1}}-  \sin\phi \sin\lambda \mathbf{i_{2}}+ \cos\phi \mathbf{i_{3}}.
\end{aligned}
\ee
The unit normal vector to the sphere at the same point $x \in \mathbb{S}^{2}$ is given by
\bel{eq:4.2}
n(x)=\cos\phi \cos\lambda \mathbf{i_{1}}+\cos\phi \sin\lambda  \mathbf{i_{2}}+\sin\phi \mathbf{i_{3}}.
\ee

In spherical coordinates, for any vector field $F$ represented by $F:= F_{\lambda } \mathbf{i_{\lambda}}+F_{ \phi }\mathbf{i_{\phi}}$, the equation of conservation law ($\ref{Eq2.4}$), can be rewritten as

\bel{eq:4.3}
 \del_t u+\frac{1}{\cos \phi } (\frac{\partial  }{\partial \phi}(F_{\phi } \cos \phi) +\frac{\partial F_{\lambda }}{\partial \lambda }   )=0.
\ee

The three general structures of the computational cells used as part of the discretization grid on the sphere are shown in Figure $\ref{Fig3}$ (a-b-c). When we go from the equator to the north or south poles, the cells are changed by a ratio of 2 at some special latitude circles to reduce the number of cells in order to satisfy the stability condition and to ensure consistency of precision in the entire domain of the sphere.

The domain of each cell $\Omega$ is defined as $\Omega : = \left \{  (\lambda, \phi) ,  \   \lambda_{1} \leqslant \lambda \leqslant \lambda_{2}, \phi_{1} \leqslant \phi \leqslant \phi_{2} \right \}$. Near the north or south poles, a triangular cell is considered which is 
a special case of the standard rectangular cell shown in Figure $\ref{Fig3}$ with zero length for the side located on the pole.

\subsection{A non-oscillatory piecewise-linear reconstruction}
\label{sec4.2}
\begin{figure}[htbp]
\begin{centering}
\includegraphics[width=8cm,height=5cm,angle=0]{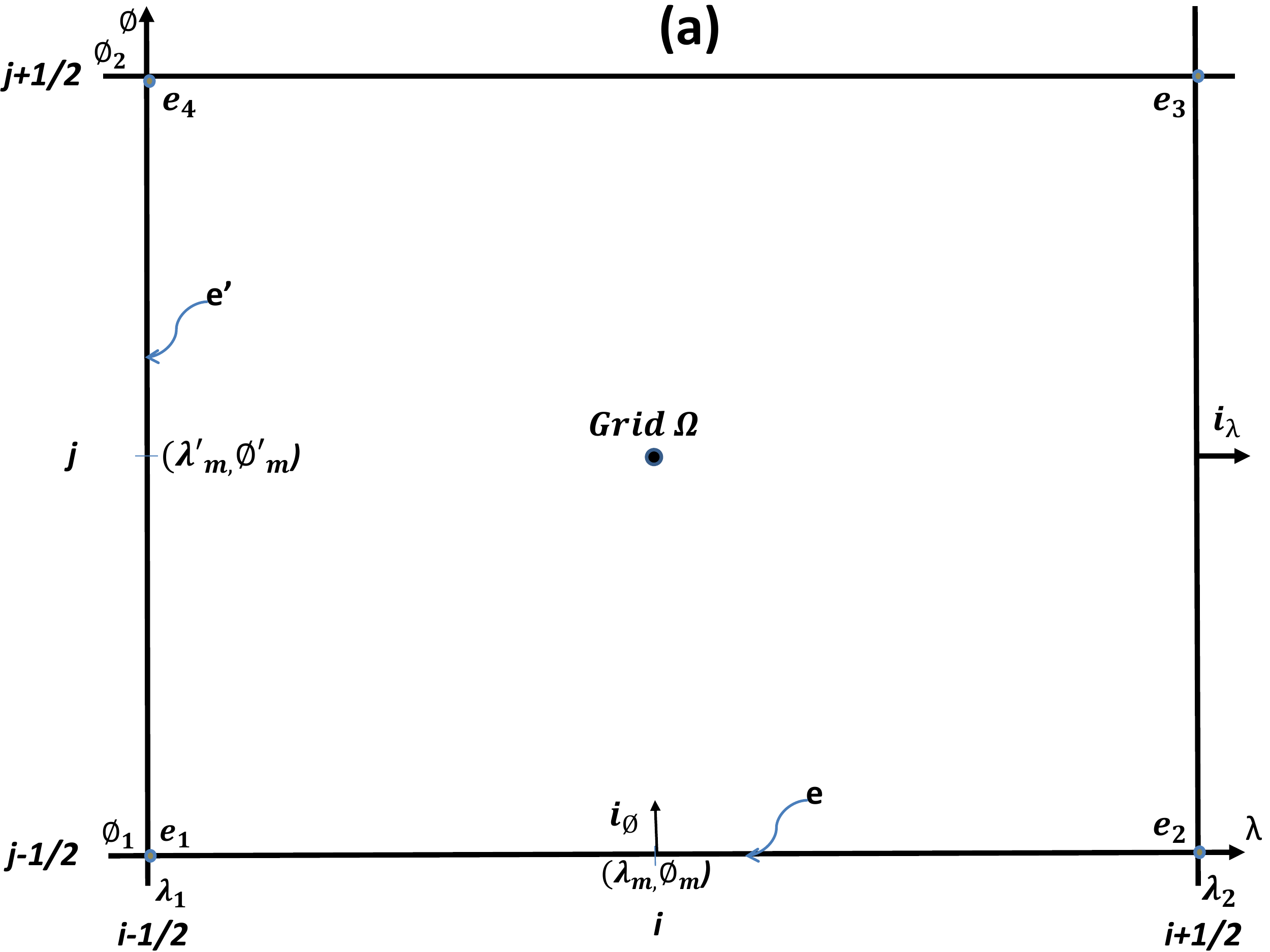} 
\includegraphics[width=8cm,height=5cm,angle=0]{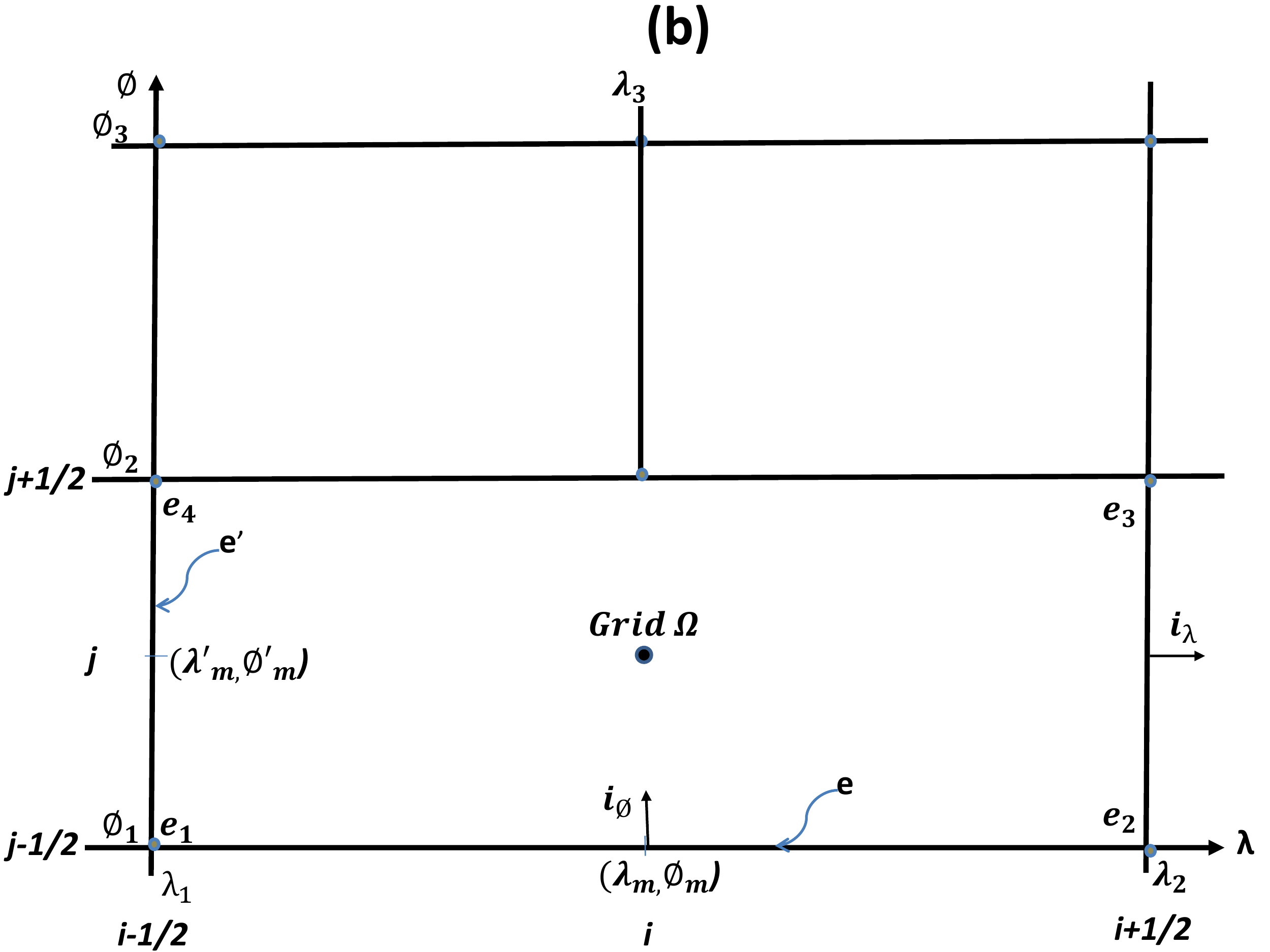} 
\includegraphics[width=8cm,height=5cm,angle=0]{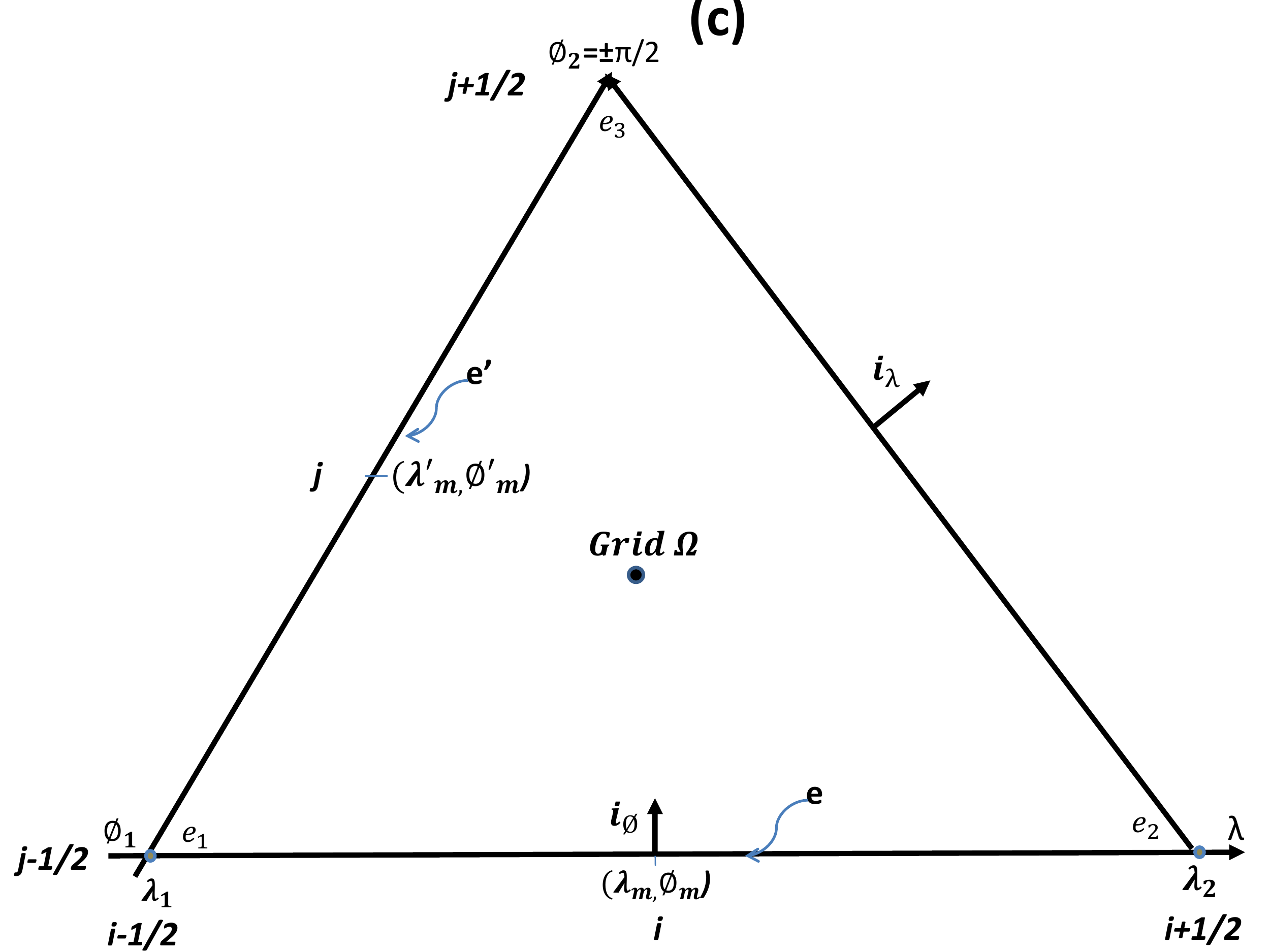} 
\includegraphics[width=8cm,height=6cm,angle=0]{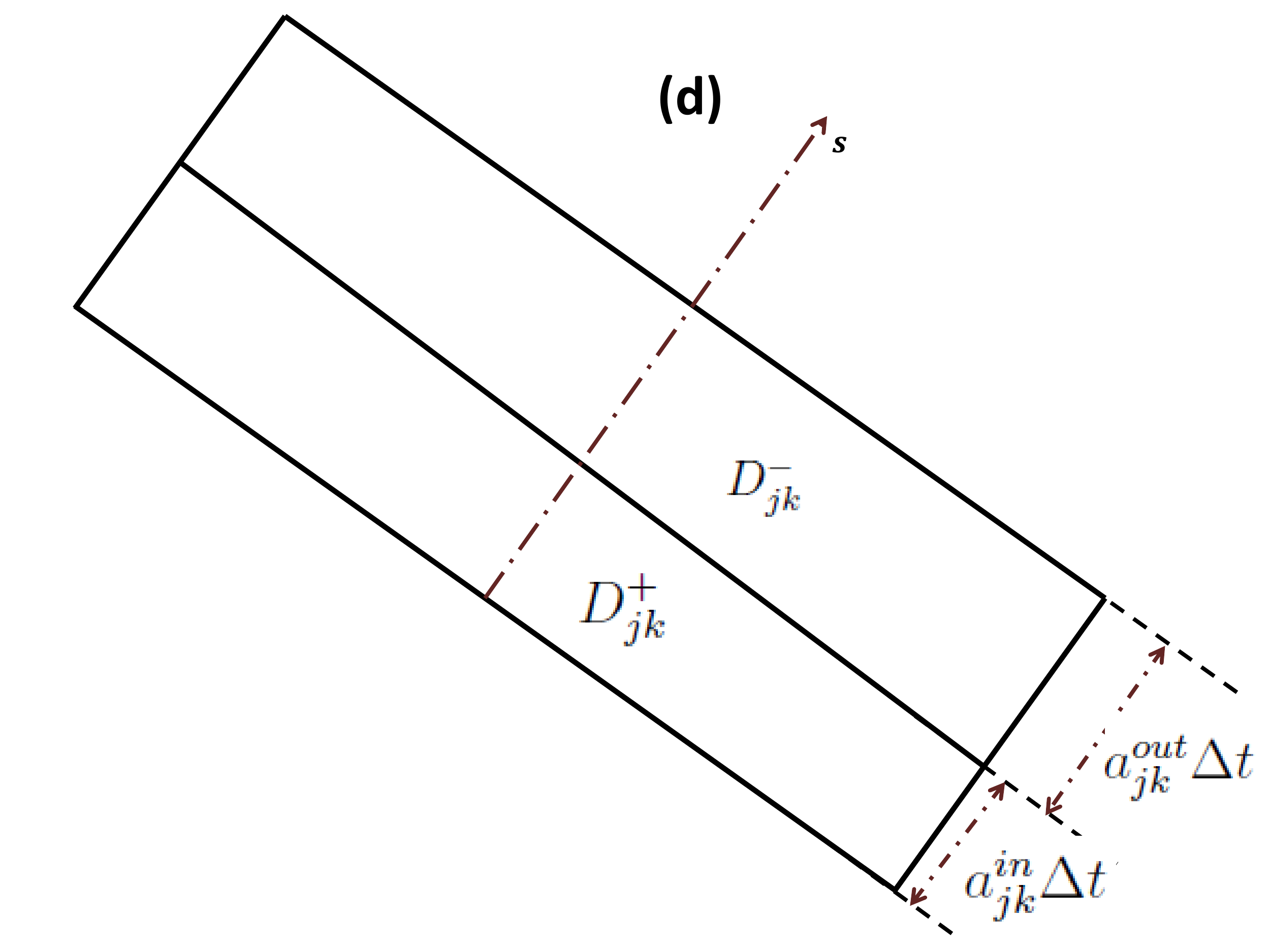} 
\par\end{centering}
\caption{ (a)-(b)-(c): Types of grids used on the sphere. (d): The domain $D_{jk}=D_{jk}^{+}\cup D_{jk}^{-}$.}
\label{Fig3} 
\end{figure}

In this section, we describe the piecewise linear reconstruction used in the proposed scheme. For simplicity, in the notations we will use the indices  $i$ and $j$ for the cell centers along the longitude and latitude, respectively (see, Figure \ref{Fig3}). 
At each time step $t_{n}$, data cell average values $u_{i,j}^{n}$ in each cell of center $(\lambda_{i},\phi_{j})$ are locally replaced by a piecewise linear function. The obtained reconstruction is as follows
\bel{eq:4.5}
u_{i,j}^{n}(\lambda,\phi) = u_{i,j}^{n} +(\lambda -\lambda_{i})\mu^{n}_{i,j}+(\phi -\phi_{j}) \sigma  ^{n}_{i,j},
\ee
where $\mu^{n}_{i,j}$ and $\sigma^{n}_{i,j}$ are the slopes in the directions of longitude and latitude, respectively. To prevent oscillations, we propose the following minmod-type reconstruction to obtain the slopes in the longitude and latitude directions
\bel{eq:4.6}
\begin{aligned}
&\mu^{n}_{i,j}=minmod [\frac{u_{i+1,j}^{n}-u_{i,j}^{n}  }{\lambda_{i+1}-\lambda_{i}}, \frac{u_{i+1,j}^{n}-u_{i-1,j}^{n}  }{\lambda_{i+1}-\lambda_{i-1}}, \frac{u_{i,j}^{n}-u_{i-1,j}^{n}  }{ \lambda_{i}-\lambda_{i-1}}  )],\\
&\sigma  ^{n}_{i,j}=minmod[\frac{u_{i,j+1}^{n}-u_{i,j}^{n} }{\phi_{j+1}-\phi_{j}}, \frac{u_{i,j+1}^{n}-u_{i,j-1}^{n} }{\phi_{j+1}-\phi_{j-1}}, \frac{u_{i,j}^{n}-u_{i,j-1}^{n} }{\phi_{j}-\phi_{j-1}}  )],
\end{aligned}
\ee
where the $minmod$ function is defined as 
\bel{eq:4.7}
\begin{aligned}
& minmod(\kappa  _{1},\kappa  _{2},\kappa  _{3})
\\
& =\begin{cases}
\kappa   \min( |\kappa  _{1}|, |\tau _{2}|,|\kappa  _{3}|), \    if   \     \kappa   =sign(\kappa  _{1})= sign(\kappa  _{2})=sign(\kappa  _{3}),   \\ 
0,  \hspace{3.5cm}   \text{otherwise.}
\end{cases}.
\end{aligned}
\ee

At each step, we compute the average values of the state variable $u$ in the computational cells. The same values are used as the values of $u$ at the cell centers of coordinates $(\lambda_{i},\phi_{j})$. The suitable points, inside the cells which 
respect these conditions for the linear reconstruction used in this study, should have the following spherical coordinates
\bel{eq:4.8}
\begin{aligned}
&\lambda _{i}=\frac{\lambda _{1}+\lambda_{2}}{2},
\\
& \phi _{j}=\frac{\phi _{2}\sin(\phi _{2})-\phi _{1}\sin(\phi _{1})+\cos \phi _{2}-\cos \phi _{1}}{\sin \phi _{2}-\sin \phi _{1}},
\end{aligned}
\ee
where $\lambda _{1}$, $\lambda _{2}$, $\phi _{1}$, and $\phi _{2}$ correspond to the longitude and latitude coordinates of the cell nodes as shown in Figure \ref{Fig3}.


\section{Geometry-compatible flux vectors and particular solutions of interest}
\label{sec5}

\subsection{Classes of geometry-compatible flux vectors}
\label{sec5.1}
 We have introduced in \cite{BLM}, two classes of flux vector fields for Equation $(\ref{eq:4.3})$. In this classification, the structure of the potential function $h(x,u)$ was used to distinguish between \textit{foliated} and \textit{generic} fluxes. In the proposed classification the parameterized level sets defined by $\Gamma _{C, u} = \big\{x\in \RR^3 \, \big/ \, h(x,u)=C \big\}$, where 
$C\in \mathbb{R}$, are used for the flux vector $F(x,u) = n(x) \wedge \nabla h(x,u)$ associated to the potential function $h$. The flux $F$ is called a \textit{foliated flux field}
 if the associated family of level sets $\big\{ \Gamma _{C, u} \big\}_{C \in \RR}$ in $\RR^3$ is independent of the state variable $u$. In other words, for any two parameters $u_1, u_2$ one can find two real numbers $C_1, C_2$ such that $\Gamma _{C_1, u_1} =  \Gamma _{C_2, u_2}$. For the generic flux field, the potential function $h=h(x,u)$ does not have this structure.

The dependency of the potential function on the space variable $x$ generates the propagation of the waves, while the dependency on the state variable $u$ leads to the formation 
of shocks in the solutions. The foliated flux with linear behavior generates the spatially periodic solutions while the foliated flux with nonlinear behavior can generate nontrivial stationary solutions.
In our analysis in \cite{BLM}, we have concluded that the new classification introduced and the character of linearity of the flux are sufficient to predict the late-time asymptotic behavior of the solutions. 
For a linear foliated flux, the solutions are simply transported along the level sets. The generic flux generates large variations in solutions, which converge to constant values within independent domains on the sphere. 
For the nonlinear foliated flux, the solution converges to its constant average in each level set. For this flux, any steady state solution should be constant along each level set. This type of 
nontrivial stationary solutions are used in our numerical experiments to demonstrate the performance of the proposed central-upwind finite volume method. 

\subsection{Particular solutions of interest}
\label{sec5.2}

The non-trivial steady state solutions which will be used in our numerical experiments are obtained using nonlinear foliated fluxes. We are particularly interested in nonlinear foliated fluxes 
based on a scalar potential function of the form 
\bel{eq:5.1}
h(x, u) =\varphi(x \cdot a)f (u), 
\ee
where $x \cdot a$ denotes the scalar product of the vector $x$ and some constant vector $a=(a_{1}, a_{2}, a_{3})^{T}  \in \mathbb{R}^{3}$, while $f$  is a function of the state variable $u$ and $\varphi $ is a function of one variable.  
This scalar potential function leads to the gradient-type flux vector field $\Phi (x,u)=\varphi ^\prime (x \cdot a) f(u) a$, where $\varphi^\prime$ is the derivative of the function $\varphi$. The flux is obtained using Equation (\ref{Eq2.2}) as
\bel{eq:5.1}
F(x,u)=\varphi ^\prime (x \cdot a) f(u) n(x)\wedge a.
\ee

For this foliated flux vector and any function $\tilde{u}$ which depends on one variable, the function defined as $u_{0}(x)=\tilde{u}(x \cdot a)=\tilde{u}(a_{1}x_{1}+a_{2}x_{2}+a_{3}x_{3})$
is a steady state solution to the conservation law $(\ref{eq:4.3})$ associated to the flux vector $F(x,u)$. 
Arbitrary functions $\varphi$ and values of the vector $a$ are used to construct nonlinear foliated fluxes and the corresponding nontrivial stationary solutions. 
In the following, $\nabla$ will be used as the standard gradient operator defined using the variable $x$ and if other variables are used, they will be specified in the notation by $\nabla_{y}$ for the gradient operator using any other variable $y$. 

In order to prove that the function $u_{0}(x)$ is a steady state solution of 
$(\ref{eq:4.3})$, the Claim 3.2 in \cite{BFL} will be used. This claim states that for any smooth function $h(x,u)$ defined on $\mathbb{S}^{2}$ with the associated gradient $\Phi=\nabla h$, if the function $u_{0}$ defined on $\mathbb{S}^{2}$ 
satisfies the condition $\nabla_{y}h(y,u_{0}(x))_{\mid y=x} = \nabla H(x)$, where H is a smooth function defined in a neighborhood of $\mathbb{S}^{2}$, then the function $u_{0}$ is a steady state solution of the conservation law $(\ref{eq:4.3})$ associated to the flux vector $F(x,u)=n(x)\wedge \Phi (x,u)$. This result will be used to prove the following corollary related to nontrivial stationary solutions which are obtained using the nonlinear foliated flux vectors. 

\bcorl{cor:2.5} \textbf{(A family of steady state solutions).} 
Consider the foliated flux vector $F(x,u)=n(x)\wedge \Phi (x,u)$ with $\Phi=\nabla h$ and $h(x, u) = \varphi(x \cdot a)f (u)$, where $a=(a_{1}, a_{2}, a_{3})^{T}  \in \mathbb{R}^{3}$, $f$  is a function of the state variable $u$ and the function $\varphi$ 
depends on one variable. For any function $\tilde{u}$ which depends on one variable, the function defined as $u_{0}(x)=\tilde{u}(x \cdot a)=\tilde{u}(a_{1}x_{1}+a_{2}x_{2}+a_{3}x_{3})$
is a steady state solution to the conservation law $(\ref{eq:4.3})$ associated to the flux $F(x,u)$. 
\ecor


We prove this result as follows, and consider the function
\bel{eq:5.3}
H(x)=H_{0} (a_{1}x_{1}+a_{2}x_{2}+a_{3}x_{3}), 
\ee
where $H_{0}$ is defined by
\bel{eq:5.4}
H_{0}(\mu )=\int_{\mu _{0}}^{\mu } \varphi^\prime (\mu ) f(\tilde{u}(\mu ))d\mu,
\ee
for some reference value $\mu _{0}$. 

The function $h(x, u) = \varphi(x \cdot a)f (u)$ is smooth in $\mathbb{R}^{3}$ and one obtains 
\bel{eq:5.5}
\nabla_{y}h(y,u_{0}(x))_{\mid y=x} =   \varphi^\prime (x\cdot a) f(\tilde{u}(x\cdot a))\sum_{k=1}^{k=3}a_{k}i_{k},
\ee
which leads to 
\bel{eq:5.6}
\nabla_{y}h(y,u_{0}(x))_{\mid y=x} = \nabla H(x)
\ee
As mentioned before, according to Claim 3.2 in \cite{BFL}, the condition (\ref{eq:5.6}) is sufficient to conclude that the function $u_{0}(x)$ is a steady state solution of the conservation law $(\ref{eq:4.3})$. This proves the claim. 

\vskip .5cm

We will consider the nonlinear foliated flux vectors based on the scalar potential functions of the form $h(x, u) =\varphi(x_{1})f (u)$, where the function $\varphi$ is not constant. For this flux, any non-constant function which depends on $x_{1}$ only is a nontrivial steady state solution of Equation $(\ref{eq:4.3})$. Another form of nonlinear foliated flux is used in our numerical tests which is obtained by using the scalar potential function of the form $h(x, u) =\varphi(x_{1}+x_{2}+x_{3})f (u)$. This case leads to steady state solutions of the form $u_{0}(x)=\tilde{u}(x _{1}+x _{2} + x _{3})$. In this paper we will consider discontinuous steady state solutions to test the performance of the proposed central-upwind method. We will use the nonlinear foliated flux vectors which are obtained by using the scalar potential function of the form $h(x, u) = \varphi(x \cdot a)f (u)$, where $f(u)=u^{2}/2$. For these flux vectors, the function defined as $u_{0}(x)=\chi (x \cdot a)\tilde{u}(x \cdot a)$ is a discontinuous stationary solution of Equation $(\ref{eq:4.3})$, where $\chi (x \cdot a)=\pm1$. 

In Tests 7 and 8, the proposed central-upwind scheme is employed to compute \textit{confined solutions} of the conservation law $(\ref{eq:4.3})$. In these cases, we consider the flux vector $F(x,u)$ which vanishes outside a domain $\Theta$ in the sphere $\mathbb{S}^{2}$. If the initial condition $u_{0}(x)$ vanishes outside of $\Theta$, then the solution should vanish outside the domain $\Theta$ for all time. However, the solution can evolve inside the domain $\Theta$ depending on the type of flux and the initial condition considered inside of $\Theta$. This case is observed in Test 7 presented in Section $\ref{sec6}$, where we choose the initial condition which is not stationary inside the domain $\Theta$ but it vanishes outside this domain. In Test 8, we will consider the flux vector $F(x,u)$ which vanishes outside $\Theta$ and is defined inside this domain using the scalar potential function $h(x, u) =\varphi(x \cdot a)f (u)$. This leads to a flux vector $F$ which satisfies the conditions mentioned in Claim \ref{cor:2.5}. For this case, we will consider an initial condition of the form $u_{0}(x)=\tilde{u}(x \cdot a)$ inside a domain $\Theta$ and vanishes outside this domain. The solution should be stationary inside $\Theta$ and should vanish outside this domain.


\section{Numerical experiments}
\label{sec6}

In this section, we demonstrate the performance of the proposed central-upwind scheme on a variety of numerical examples. Different types of nonlinear foliated fluxes are used to construct some particular and interesting solutions.
In Example 1, four numerical tests are performed using different discontinuous steady state solutions of the conservation law $(\ref{eq:4.3})$ with the nonlinear foliated flux vectors based on the scalar potential functions of the form $h(x, u) =\varphi(x_{1})f (u)$.
In Example 2, two numerical tests are performed using different discontinuous steady state solutions in the spherical cap of $(\ref{eq:4.3})$ which are obtained by using the nonlinear foliated flux corresponding to the scalar potential function of the form $h(x, u) =\varphi(x_{1}+x_{2}+x_{3})f (u)$. In the third example, two numerical tests are performed where the proposed scheme is employed to compute confined solutions. 

\subsection*{Example 1---Discontinuous steady state solutions}

First, we consider the potential function $h(x, u) =x_{1}f (u)$, where $f(u)=u^{2}/2$ which leads to the nonlinear foliated flux vector $F(x,u)=f(u) n(x)\wedge i_{1}$. We take the following 
discontinuous steady state solution of Equation $(\ref{eq:4.3})$ as initial condition (Test 1).
\bel{eq:6.1}
\begin{aligned}
& u_{2}(x)= \begin{cases}
&\gamma  x_{1}^{3},   \hspace{1.2cm} -1 \leq x_{1} \leq 0.5, \\
&- \gamma x_{1}^{2}/(2x_{1}+1) , \hspace{0.2cm}  0.5 \leq x_{1} \leq 1,
\end{cases} 
\end{aligned}
\ee
where $\gamma$ is an arbitrary constant which controls the amplitude and shocks of the solution. This solution has a single closed curve of discontinuity on the sphere.

\begin{figure}[htbp]
\begin{centering}
\includegraphics[width=8cm,height=8cm,angle=0]{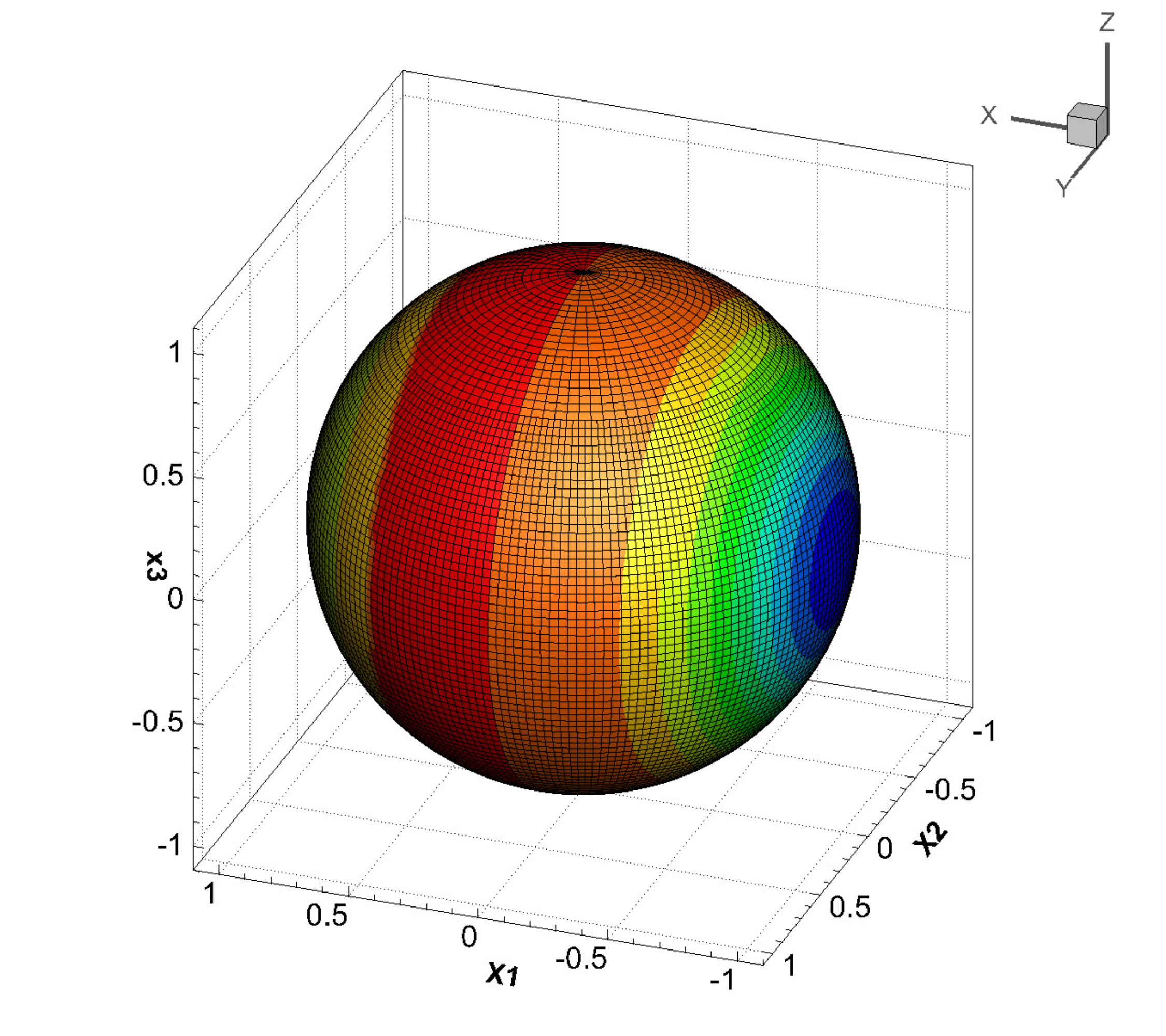} 
\includegraphics[width=8cm,height=8cm,angle=0]{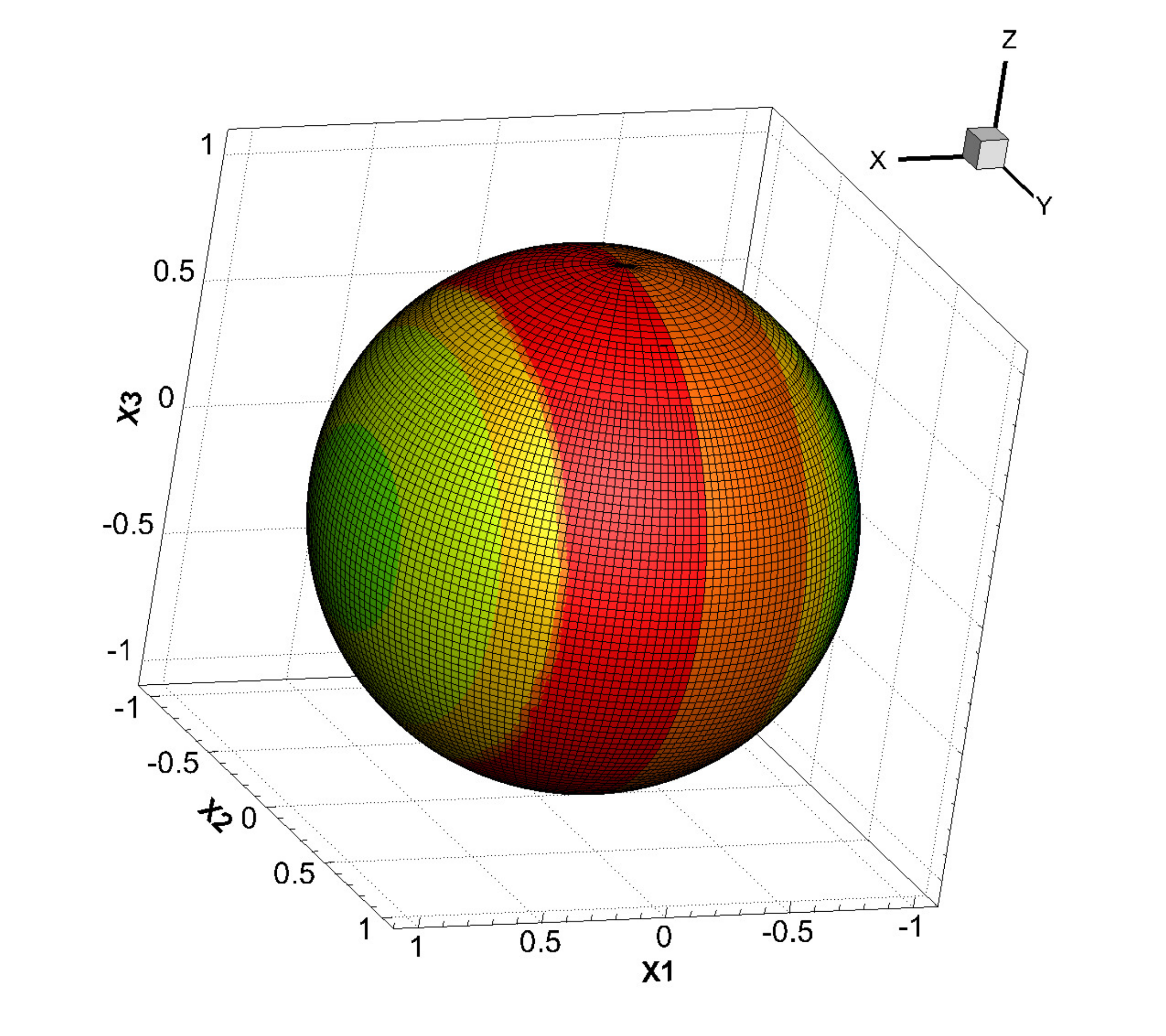} 
\par\end{centering}
\caption{Solutions on the entire sphere at time $t=5$ for Test 1 (left) and Test 2 (right) }
\label{Fig4} 
\end{figure}

The numerical solution is computed using a grid with an equatorial longitude step $ \Delta \lambda =\pi /96$ and a latitude step $\Delta\phi  =\pi /96$, and a time step $\Delta t= 0.04$. 
Figure \ref{Fig4}, on the left, shows the numerical solution with $\gamma=0.1$ which is computed using the proposed scheme at a global time $t=5$. The numerical solution remains nearly unchanged in time using the proposed scheme. The numerical solution error defined by using the $L^{2}-$norm is computed by summation over all grid cells on the sphere. For Test 1, the error is $u_{error}=1.5 \times 10^{-4}$ at time $t=5$ which is small compared to the full range of the numerical solution $u_{max}-u_{min}=0.11237$.

Another test is performed using the steady state solution (\ref{eq:6.1}) as initial condition with $\gamma=0.5$ (Test 2) and the same computational grid used in Test 1 and a time step $\Delta t= 0.04$. As is shown in Figure \ref{Fig4}, on the right for Test 2, the solution remains nearly unchanged up to a global time $t=5$. The error using the $L^{2}-$norm is $u_{error}=2.7 \times 10^{-3}$, which is small compared to the full range of the solution $u_{max}-u_{min}=0.56185$.

Now we consider a new test (Test 3) using the following steady state solution, with more discontinuities, which is defined in three domains separated by two closed curves on the sphere
\bel{eq:6.2}
\begin{aligned}
& u_{2}(x)= \begin{cases}
&\gamma  x_{1}^{4},   \hspace{1.2cm} -1 \leq x_{1} \leq -0.5, \\
& 0.5 \gamma x_{1}^{3}, \hspace{0.7cm} -0.5 < x_{1} < 0.5, \\
& -0.25\gamma  x_{1}^{2}, \hspace{0.2cm}  0.5 \leq x_{1} \leq 1.
\end{cases} 
\end{aligned}
\ee

\begin{figure}[htbp]
\begin{centering}
\includegraphics[width=8cm,height=8cm,angle=0]{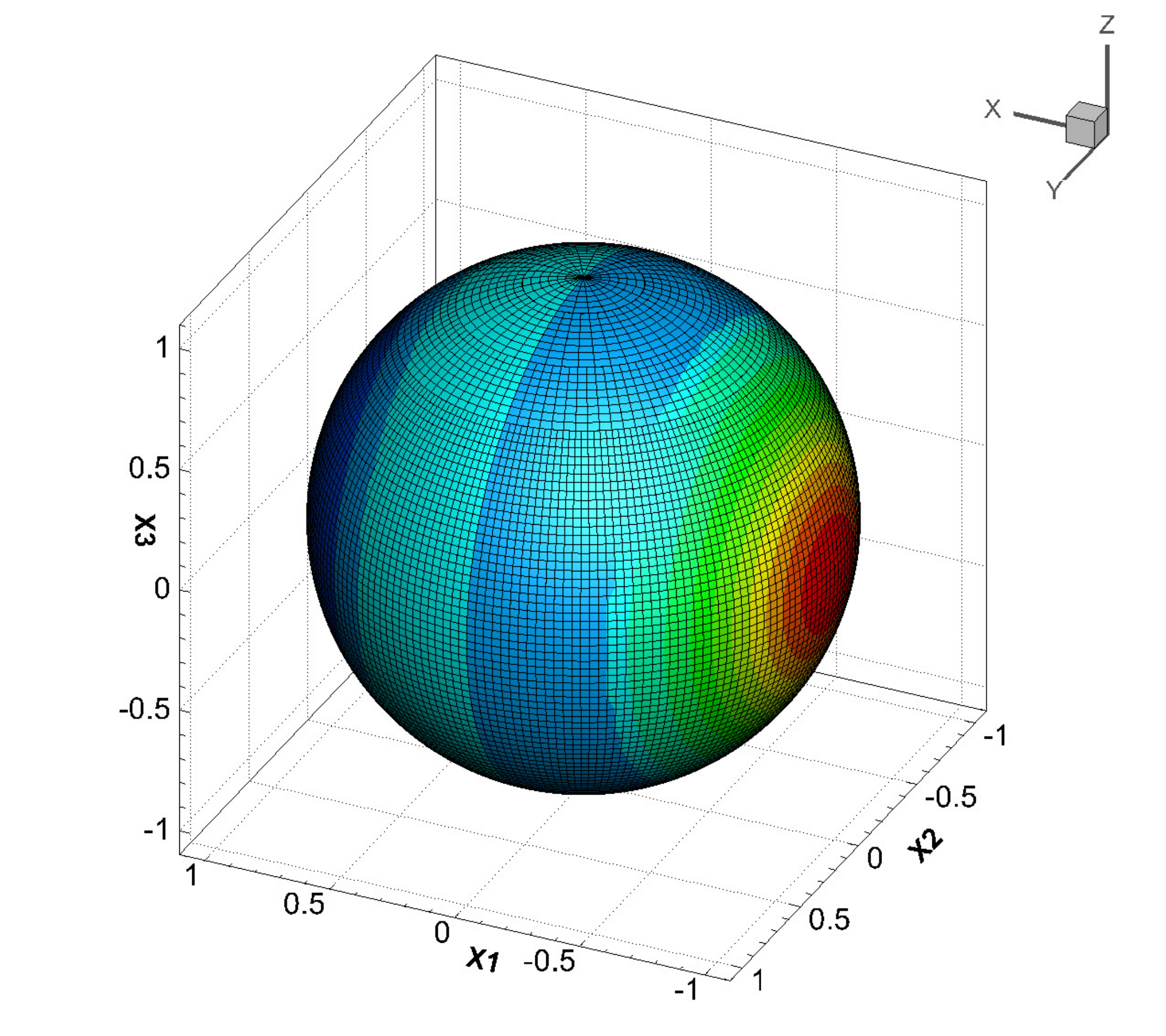} 
\includegraphics[width=8cm,height=8cm,angle=0]{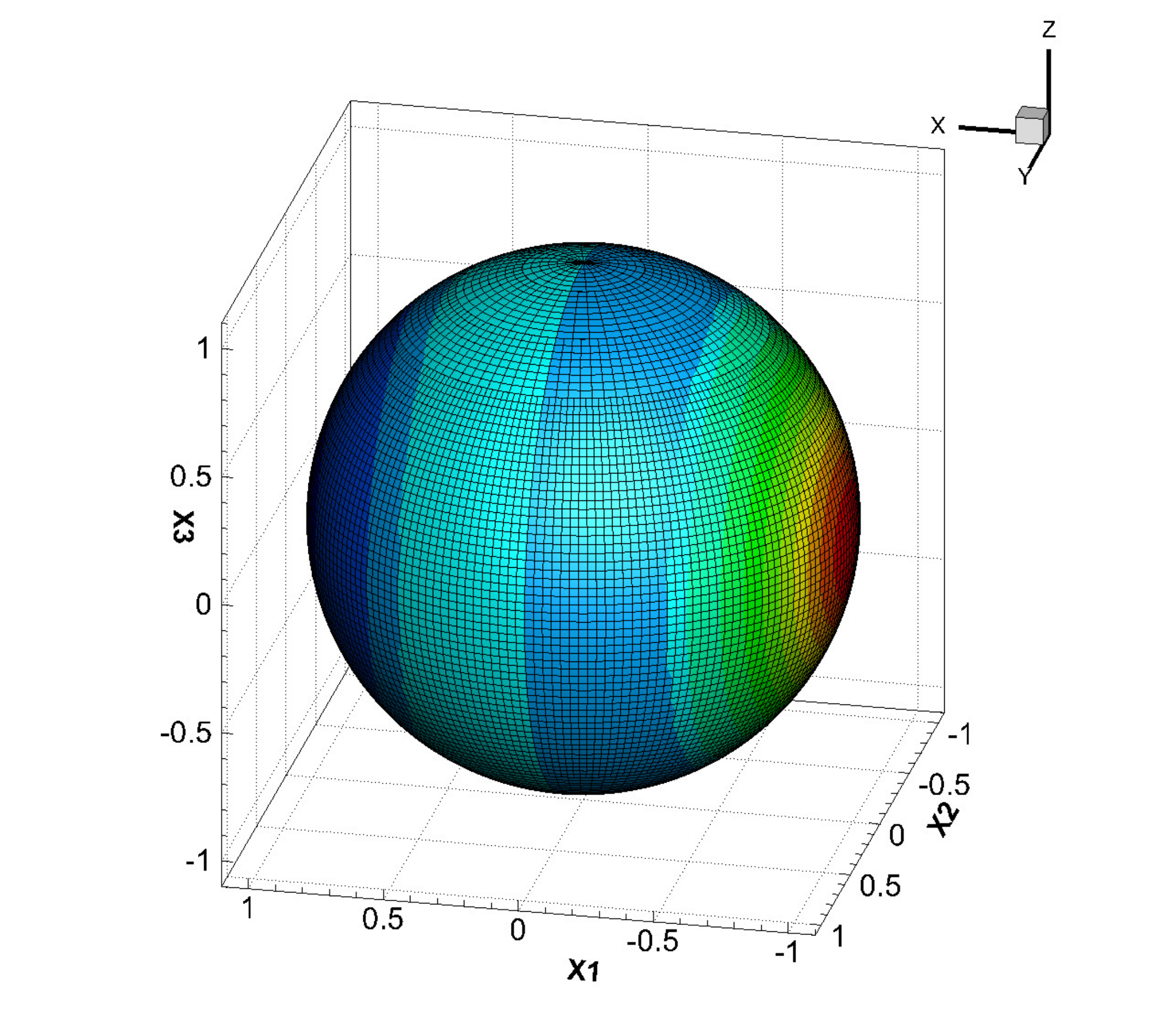} 
\par\end{centering}
\caption{Solutions on the entire sphere at time $t=5$ for Test 3 (left) and Test 4 (right) }
\label{Fig5} 
\end{figure}

The numerical solution is computed using a time step $\Delta t= 0.04$ and the same grid on the sphere used in the previous tests. As shown in Figure \ref{Fig5}, on the left, the numerical solution which is obtained at time $t=5$ using the 
proposed method based on the initial condition (\ref{eq:6.2}) with $\gamma=0.1$ remains nearly unchanged. The error is $u_{error}=9.6 \times 10^{-5}$, which is small compared to the full range $u_{max}-u_{min}=0.12488$.

For $\gamma=0.5$ (Test 4) we used the same computational grid and a time step $\Delta t= 0.04$. As is shown in Figure \ref{Fig5}, on the right, again for this test the numerical solution at time $t=5$ remains nearly unchanged. 
The error using the $L^{2}-$norm is $u_{error}=1.9 \times 10^{-3}$, which is small compared to the full range of the solution $u_{max}-u_{min}=0.62441$.


\subsection*{Example 2---Discontinuous steady state solutions in a spherical cap}

In the following, the performance of the proposed finite volume method will be analyzed using some particular steady state solutions in a spherical cap. The scalar potential function $h(x,u)= (x _{1}+x _{2} + x _{3})f(u)$ is considered with $f(u)=u^{2}/2$. This leads to the nonlinear foliated flux $F(x,u)= f(u)n(x)\wedge (i_{1}+ i_{2}+ i_{3})$.  The function of the form $u(x)=\chi (\theta) \tilde{u} (\theta)$ is a steady state solution of Equation $(\ref{eq:4.3})$, where $\tilde{u}$ is an arbitrary real function depending on one variable and 
$\theta =x_{1}+x_{2}+x_{3}$. 

In this numerical example (Test 5), the following discontinuous steady state solution is considered as initial condition
\bel{eq:6.3}
\begin{aligned}
& u(0,x)= \begin{cases}
& \text{}  0.1 / (\theta+2) ,  \hspace{1.6cm}    0 \leq \theta , \\
& \text{} -0.1 / (\theta+2) ,  \hspace{1.1cm}   \text{otherwise}. \\
\end{cases} 
\end{aligned}
\ee

\begin{figure}[htbp]
\begin{centering}
\includegraphics[width=8cm,height=8cm,angle=0]{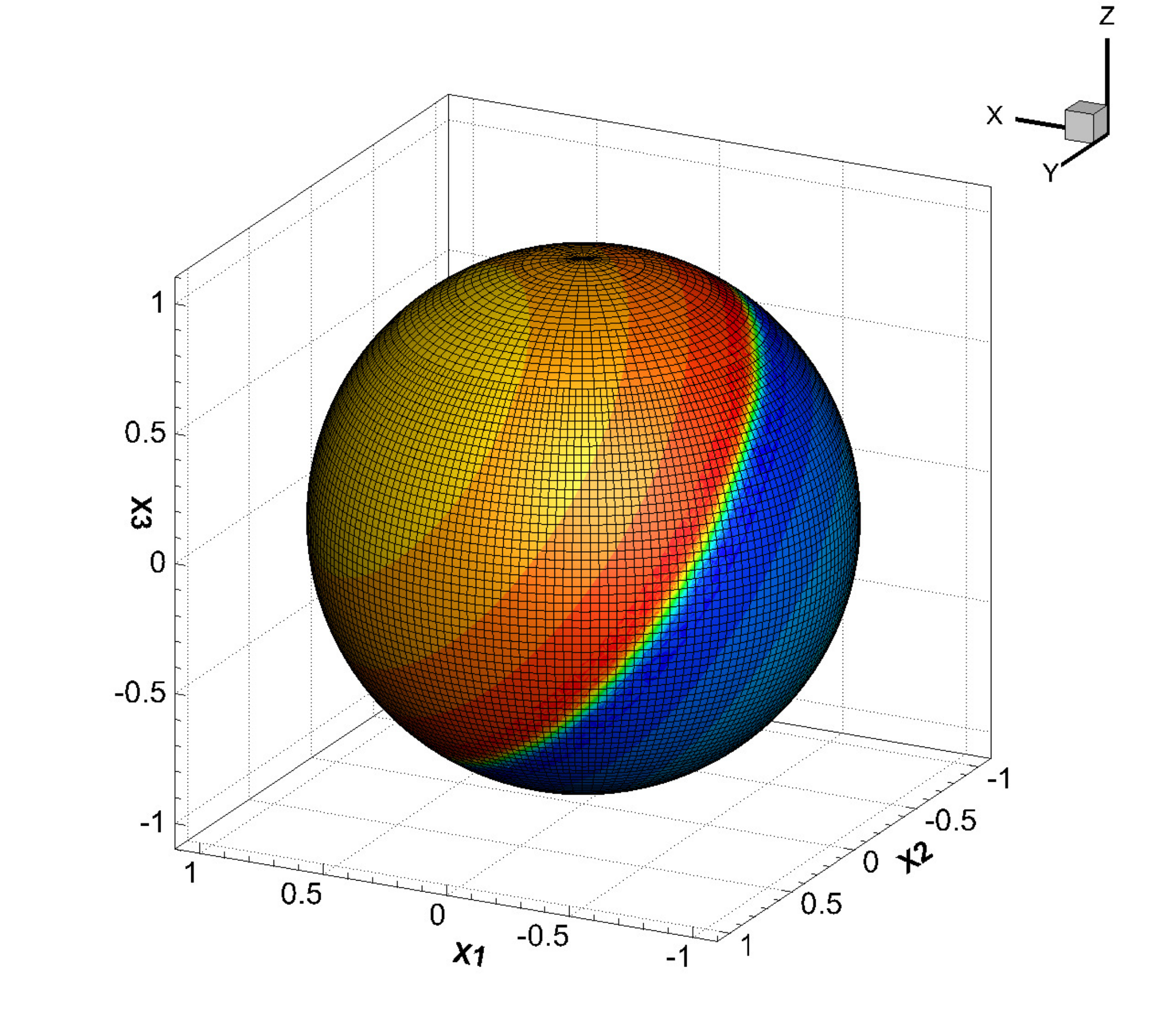} 
\includegraphics[width=8cm,height=8cm,angle=0]{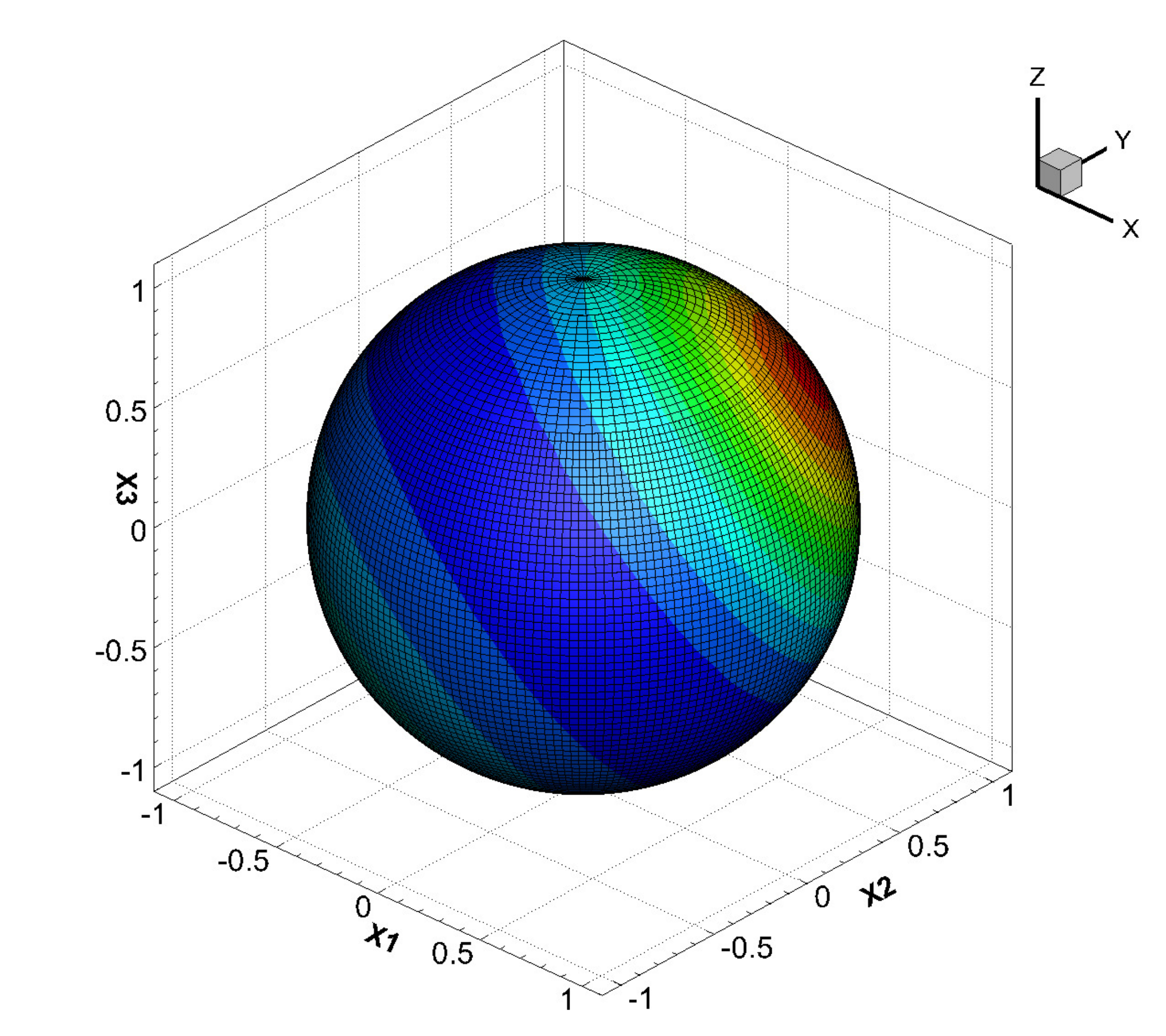} 
\par\end{centering}
\caption{Solutions on the entire sphere at time $t=5$ for Test 5 (left) and Test 6 (right) }
\label{Fig6} 
\end{figure}

The numerical solution is computed by using a grid with an equatorial longitude step $ \Delta \lambda =\pi /96$ and a latitude step $\Delta\phi  =\pi /96$, and a time step $\Delta t= 0.02$. 
Figure \ref{Fig6}, on the left, shows the numerical solution which remains nearly unchanged in time after being subjected to integration up to a global time $t=5$ by the proposed scheme. The numerical solution error defined by using the $L^{2}-$norm is 
$u_{error}=1.3 \times 10^{-3}$, which is small compared to the full range $u_{max}-u_{min}=0.4232$.

The following numerical example (Test 6) is performed using the same nonlinear foliated flux considered in Test 5 and the steady state solution with more discontinuities defined by
\bel{eq:6.4}
\begin{aligned}
& u(0,x)= \begin{cases}
& \text{}  0.2  \theta^{3} ,  \hspace{1.3cm}   0.5  \leq  \theta , \\
& \text{}  0.1  \theta^{2},  \hspace{1.3cm}    \theta  \leq  -0.5 , \\
& \text{} -0.025 ,  \hspace{0.8cm}   \text{otherwise}. \\
\end{cases} 
\end{aligned}
\ee

 The numerical solution is computed using the same grid used in Test 5 and a time step $\Delta t= 0.02$. 
Figure \ref{Fig6}, on the right, shows the numerical solution at time $t=5$ which remains stationary with the error $u_{error}=1.8 \times 10^{-3}$ which is negligible compared to the full range of the solution $u_{max}-u_{min}=1.0638$.

\subsection*{Example 3---Confined solutions}

In this part, two numerical tests are performed using confined solutions of the conservation law $(\ref{eq:4.3})$ based on the flux vector which is obtained using the following potential function

\bel{eq:6.5}
\begin{aligned}
&h(x,u)= \begin{cases}
& \text{} x_{1}^{2}f_{1}(u))  ,  \hspace{1.3cm}  x_{1} \leq 0 , \\
& \text{}0 ,  \hspace{0.8cm}   \text{otherwise}. \\
\end{cases} 
\end{aligned}
\ee 

\begin{figure}[htbp]
\begin{centering}
\includegraphics[width=8cm,height=8cm,angle=0]{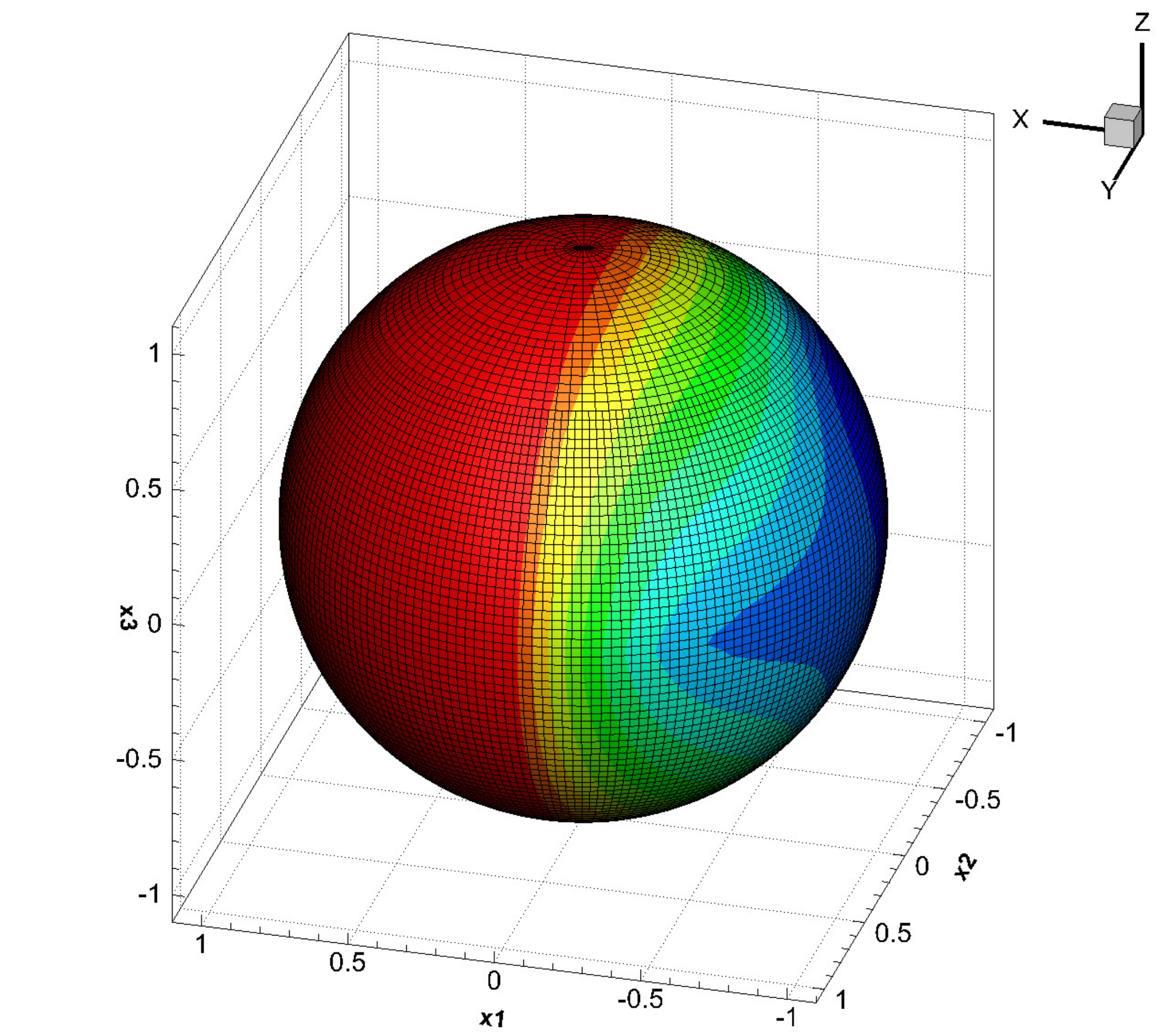} 
\includegraphics[width=8cm,height=8cm,angle=0]{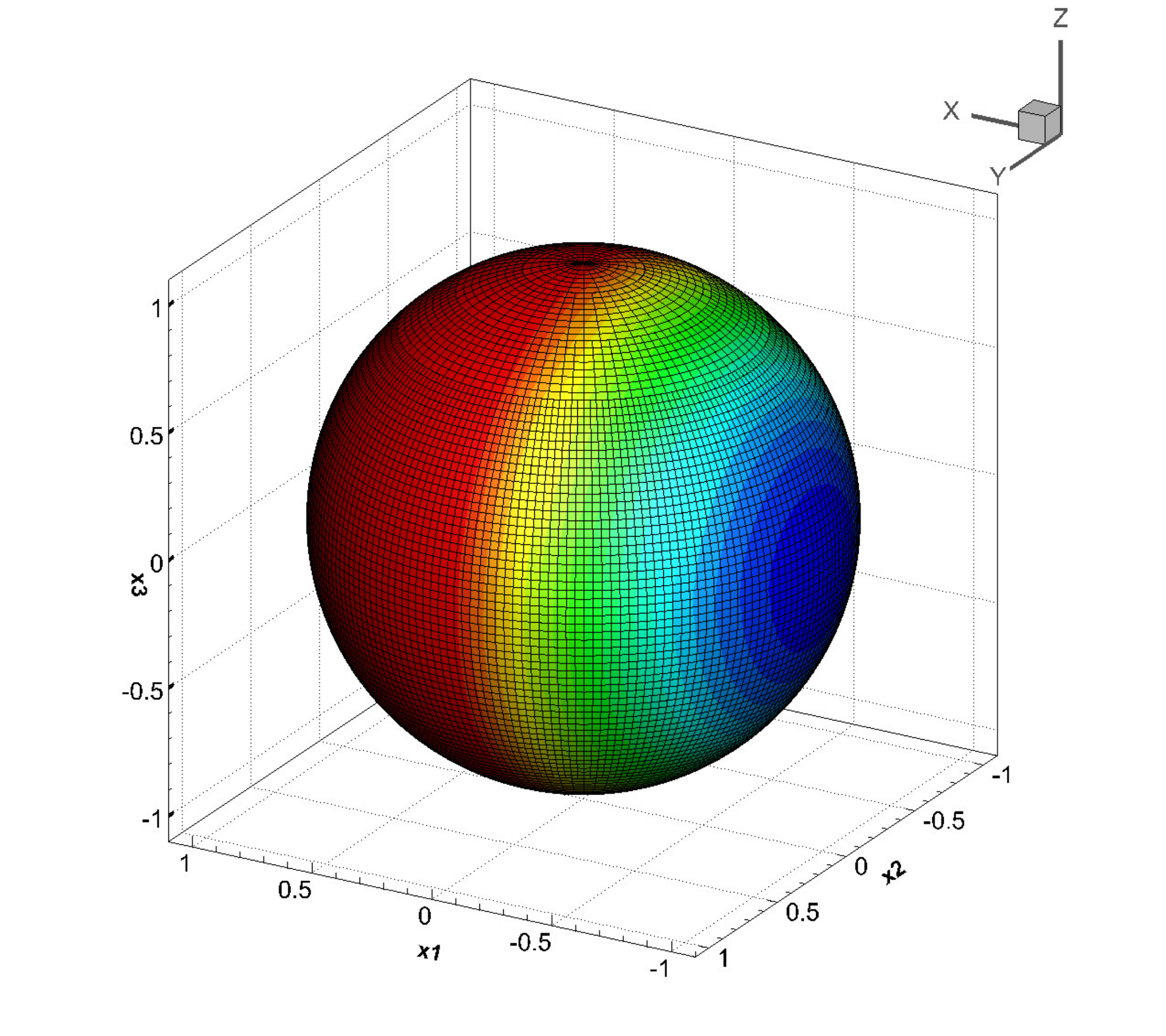} 
\par\end{centering}
\caption{Solutions on the entire sphere at time $t=5$ for Test 7 (left) and Test 8 (right) }
\label{Fig78} 
\end{figure}

In Test 7, we consider the following function
\bel{eq:6.6}
\begin{aligned}
&u(x,0)= \begin{cases}
& \text{} 0.1   (1+x_{2}^{2})x_{1} ,  \hspace{1.3cm}  x_{1} \leq 0 , \\
& \text{}0 ,  \hspace{0.8cm}   \text{otherwise}. \\
\end{cases} 
\end{aligned}
\ee
The solution of the conservation law $(\ref{eq:4.3})$, which is obtained using the function (\ref{eq:6.6}) as initial condition, is confined and it vanishes outside the domain $x_{1}\leq 0$. The numerical solution is computed using the proposed scheme 
with an equatorial longitude step $ \Delta \lambda =\pi /96$, a latitude step $\Delta\phi  =\pi /96$, and a time step $\Delta t= 0.04$. Figure \ref{Fig78}, on the left, shows the numerical solution at time $t=5$. The solution evolves in time inside the domain $x_{1} \leq 0$, but it vanishes outside this domain which is in good agreement with the evolution of the analytical solution.

In the second numerical test (Test 8), we consider an initial condition which is a confined solution and steady state inside the domain $x_{1}\leq 0$. The following initial condition is considered
\bel{eq:6.7}
\begin{aligned}
&u(x,0)= \begin{cases}
& \text{} 0.1 x_{1} ,  \hspace{1.3cm}  x_{1} \leq 0 , \\
& \text{}0 ,  \hspace{1.0cm}   \text{otherwise}. \\
\end{cases} 
\end{aligned}
\ee
The numerical solution is computed using the proposed central-upwind scheme 
with the same grid and time step which are used in Test 7. Figure \ref{Fig78}, on the right, shows the numerical solution at time $t=5$. The solution remains steady state in the domain $x_{1} \leq 0$, and it vanishes outside this domain for all time as the initial condition which is in good agreement with the evolution of the analytical solution. The $L^{2}$-error of the numerical solution over the sphere is $u_{error}=9.6 \times 10^{-5}$ at time $t=5$, which is small compared to the full range of the solution $u_{max}-u_{min}=0.2$


\section{Concluding remarks}
\label{sec7}

We have introduced a new geometry-preserving, central-upwind scheme for the discretization of hyperbolic conservation laws posed on the sphere.
The main advantage of the proposed scheme is its simplicity since it does not use Riemann problems. The proposed scheme is strongly connected 
to the analytic properties of the equation and the geometry of the sphere. In the proposed scheme, in order to improve the accuracy, the Gaussian quadrature can be used instead of the midpoint rule to compute the spatial integrals. The Gaussian quadrature will not have any impact on the geometry-compatibility condition of the scheme. 
The semi-discrete form of the proposed method using the Gaussian quadrature will remain strongly connected 
to the analytic properties of the equation and the geometry of the sphere. 

In the proposed method, a non-oscillatory
reconstruction is used in which the gradient of each variable is computed using  a minmod-function to ensure stability. Our numerical experiments demonstrate the ability of the proposed scheme to avoid
oscillations. The performance of the second-order version of the designed scheme is tested using numerical examples. The results clearly demonstrated the proposed scheme's potential and ability to 
 resolve the discontinuous solutions of conservation laws on the sphere. 

Note that the formulation of the semi-discrete form of the proposed method is based on some approximations and assumptions. The scheme is more suitable for discontinuous solutions with shocks 
of average amplitude. However, the proposed method has the advantage of simplicity compared to upwind schemes. As previously mentioned, the first advantage is that the proposed scheme is Riemann-problem-solver-free. The second advantage is related to the 
resolution, where the proposed scheme does not use any splitting approach which is widely used in upwind schemes to simplify the resolution of the Riemann problem. This again renders the proposed numerical scheme less expensive compared to 
upwind methods. The scheme developed for scalar nonlinear hyperbolic conservation laws could be extended to multidimensional hyperbolic conservation laws and shallow water models posed on the sphere.



\end{document}